\documentclass[11pt]{article}

\usepackage{sec}
\usepackage{makeidx}
\usepackage{diagram}
\usepackage{psfig}
\usepackage{amssymb}

\newtheorem{theorem}{Theorem}[section]
\newtheorem{lemma}[theorem]{Lemma}
\newtheorem{proposition}[theorem]{Proposition}
\newtheorem{corollary}[theorem]{Corollary}
\newtheorem{remark}[theorem]{Remark}

\newtheorem{ex}[theorem]{Example}

\newcommand\bC{{\mathbb C}}

\newcommand{\bQ}{{{\mathbb Q}}}
\newcommand\bR{{\mathbb R}}
\newcommand\bS{{\mathbb S}}

\newcommand\bZ{{\mathbb Z}}

\newcommand{\dir}{{\mathfrak D}}

\newcommand{\ii}{{\bf i}}
\newcommand{\bc}{{\bf c}}

\newcommand{\p}{{\cal P}}

\newcommand{\co}{{\sf C}}

\def\ra{\rightarrow}

\def\be{\begin{equation}}
\def\ee{\end{equation}}

\def\lan{\langle}
\def\ran{\rangle}
\newcommand{\LL}{\left(\left( }
\newcommand{\RR}{\right)\right)}

\renewcommand{\Box}{\blacksquare}

\newcommand{\si}{\sigma}

\newcommand{\vfi}{{\varphi}}

\newcommand{\mod}{ {{\rm mod}\,}}

\newcommand{\one}{{\bf{\hat{1}}}}

\begin{document}

\title{Seiberg-Witten Theoretic Invariants of Lens Spaces  }

\author{Liviu I. Nicolaescu\\Department of Mathematics\\University of Notre Dame\\Notre Dame, IN 46556\\nicolaescu.1@nd.edu}

\date{Version 3-January 5, 1999}



\maketitle
\begin{abstract} We  describe an effective algorithm for computing  Seiberg-Witten
invariants of lens spaces. We  apply it to two problems: (i) to
compute the Froyshov invariants of a large  family of  lens
spaces; (ii)  to show that the knowledge of the Seiberg-Witten
invariants of a lens space  is topologically equivalent to the
knowledge of  its Casson-Walker invariant  and of its
Milnor-Turaev torsion.

 Problem (i) has several interesting topological consequences concerning the
negative definite manifolds  bounding a given lens space.
\end{abstract}



\addcontentsline{toc}{section}{Introduction}

\begin{center}
{\bf Introduction}
\end{center}

\bigskip

The Seiberg-Witten theory of rational homology spheres is
particularly difficult since  the usual count of monopoles  leads
to a metric {\em dependent} integer.   W. Chen \cite{Chen1}, Y.
Lim \cite{Lim} and M. Marcolli \cite{Mar} have   shown that this
count,  suitably altered  by  a certain combination of eta
invariants,  leads to a topological invariant.     For integral
homology spheres, there is an unique $spin^c$ structure and this
altered count was shown to coincide with the Casson invariant; see
\cite{Chen2}, \cite{Lim2},  and \cite{N2} in the special case of
Brieskorn spheres.     For a rational homology sphere $N$ there
are $\# H_1(N, {\bZ})$ such invariants which are rational numbers.
They define a function
\[
{\bf sw}: Spin^c(N)\ra {\bQ},\;\; \si \mapsto {\bf sw}(\si).
\]
We will call ${\bf sw}$ the Seiberg-Witten invariant of $N$. This invariant can be further formalized as follows.

Recall that $H_1(N, {\bZ})\cong H^2(N, {\bZ})$ acts freely and
transitively on the space  $Spin^c(N)$ of $spin^c$ structures on
$N$
\[
 Spin^c(N) \times H_1(N, {\bZ})  \ni (\si,h)\mapsto \si\cdot h \in Spin^c(N)
\]
Thus each $\si_0\in Spin^c(N)$ defines an
element ${\bf SW}_{\si_0}\in {\bQ}[H]$ ($=$ the rational group algebra
of the multiplicative group $H=H_1(N, {\bZ})$) defined by
\[
{\bf SW}_{\si_0}=\sum_{h\in H}{\bf sw}(\si_0\cdot h)h.
\]
Clearly
\[
{\bf SW}_{\si_0\cdot g}= {\bf SW}_{\si_0}\cdot g^{-1},\;\;\forall
g\in H.
\]
Thus, the   collection ${\bf SW}:=\{{\bf SW}_\si;\;\;\si \in Spin^c(N)\}\subset
{\bQ}[H]$  coincides with an orbit of the   right action of $H$ on
${\bQ}[H]$ so that the Seiberg-Witten invariant can be viewed  as
an element in ${\bQ}[H]/H$.

This Seiberg-Witten invariant is unchanged by natural involution $\bar{}:Spin^c(N)\ra Spin^c(N)$, $\si \mapsto
\bar{\si}$. The conditions ${\bf sw}(\si)={\bf sw}(\bar{\si})$  and
$\bar{\si\cdot h}=\bar{\si}\cdot h^{-1}$ imply
\[
{\bf SW}_{\bar{\si}}=\overline{SW_\si}
\]
where $\bar{}:{\bQ}[H]\ra {\bQ}[H]$  is the  involution
determined by $H\ni h\mapsto h^{-1}\in H$.

A few years ago, using Seiberg-Witten theory,    Kim Froyshov (\cite{Fr}) defined another invariant
of a rational homology $3$-sphere $N$ which contains nontrivial information about the possible negative definite  $4$-manifolds
which can bound $N$.    His invariant  is the sum of a highly
unorthodox count of solutions of the Seiberg-Witten equations  and   the same combination of eta invariants   entering
the definition of ${\bf sw}$.    This is done for each $spin^c$
structure and then, the maximum amongst these numbers is chosen.
The resulting  rational number is still metric dependent. To get rid of
this dependence  Froyshov takes  the infimum over all
``reasonable'' metrics on $N$.

In \cite{N2} we have  explicitly computed the invariant ${\bf SW}$
for Brieskorn homology spheres with at most $4$ singular fibers
and we have  identified it with the Casson invariant. In \cite{N}
we  computed Froyshov's invariant  for many   Brieskorn homology
spheres with $3$ singular fibers and we have indicated an
algorithm for producing upper estimates  for any Brieskorn sphere
with three singular fibers.

In the present paper we use the results and techniques of \cite{N}
to produce   a simple algorithm computing  the  ${\bf SW}$ and the
Froyshov invariants of lens spaces.   As in \cite{N}, these
formul{\ae} involve the Dedekind-Rademacher sums so,  each
concrete computation, although completely elementary,   can be
quite involved.  On the positive side, these computations can be
easily implemented on any computer algebra system  (such as {\it
MAPLE})  and the numerical  experiments reveal  very beautiful
patterns (see (\ref{eq: liviu0}), (\ref{eq: liviu2})-(\ref{eq:
liviu9})) and \S 3.3. The concrete computations lead to
interesting topological consequences. Here are  some samples of
them.

\medskip

{\em The lens space $L(2k+1,1)$ bounds no smooth, even, negative definite $4$-manifold while the lens space $L(4k+1,2)$ bounds no smooth, even, negative definite $4$-manifold $X$ such that $H_1(X,{\bZ})$ has no $2$-torsion.}

\medskip

Using recent results of Paolo Lisca, \cite{Lisca},  we can    deduce some information about the fillable contact structures on lens spaces.  The following is an immediate consequence of Lisca's work and our computations.

\medskip

{\em $L(2k+1,1)$   cannot be the contact boundary of any even,
symplectic  $4$-manifold.}

\medskip

Denote  by ${\bf SW}_{p,q}$ is the Seiberg-Witten  invariant of
$L(p,q)$.  It is an element  of ${\bQ}[{\bZ}_p]/{\bZ}_p$ and we
will regard it as a polynomial in one variable $t$ satisfying
$t^p=1$.

The ring ${\bQ}[{\bZ}_p]$ is equipped with an augmentation map
\[
{\bf aug}: {\bQ}[{\bZ}_p]\ra
{\bQ},\;\;\sum_{k=0}^{p-1}a_kt^k\mapsto\sum_{k=0}^{p-1}a_k.
\]
We prove in  \S 3.2, Theorem  \ref{th: cw}  that
\be
{\bf aug}({\bf SW}_{p,q})= CW(L(p,q)).
\label{eq: cw}
\ee
{\em where $CW$ denotes the  Casson-Walker
invariant (see \cite{Walker})  of a rational homology sphere normalized as in \cite{Lescop}.}

Following \cite{Milnor}  we introduce the polynomial
\[
\Sigma=\sum_{k=0}^{p-1}t^k.
\]
It can be used to define a projection
\[
{\bf Proj}:{\bQ}[{\bZ}_{p}]\ra \Lambda_p :=\ker {\bf aug},\;\;  R\mapsto R-\frac{{\bf aug}(R)}{p}\Sigma.
\]
Set
\[
T_{p,q}={\bf Proj}({\bf SW}_{p,q})= {\bf
SW}_{p,q}-\frac{CW(L(p,q))}{p}\Sigma.
\]
We can regard $T_{p,q}$ as an element of $\Lambda/{\bZ}_p$. If
$A$, $B$ are two  ``polynomials'' in $\Lambda_p$ then $A\sim B$
will signify $A=t^nB$ for some $n\in {\bZ}$.

The Milnor torsion  of $L(p,q)$,  which we denote by $\tau_{p,q}$,  is also an  element  of
$\Lambda_p$  (see \cite{Milnor}). More precisely,  using the
convention of \cite{Turaev} we have (see \cite{Milnor},
\cite{Turaev})
\[
\tau_{p,q} \sim (1-t)^{-1}(1-t^q)^{-1}
\]
i.e.
\[
\tau_{p,q}(1-t)(1-t^q)\sim {\one}:=1-\frac{1}{p}\Sigma.
\]
As explained in \cite{Milnor} the ``polynomial''  $\one$
represents {\bf 1} in $\Lambda_p$. We prove the following.

\medskip

{\em For any   positive integers $p,q$ such that
$g.c.d(p,q)=g.c.d.(p,q-1)=1$ we have}
\be
T_{p,q}(1-t)(1-t^q)\sim \tau_{p,q} \label{eq: turaev} \ee

\medskip

The restriction $g.c.d.(p,q-1)=1$ is  purely for technical
reasons, to slightly simplify certain accounting jobs. The method
we present works in the general case, when $g.c.d(p,q-1)\geq 1$ .
We did not consider the details to be very enlightening so we have
not included them. We will present them elsewhere. Instead,bwe
present the results of some numerical experiments confirming the
equality $T_{p,q}\sim \tau_{p,q}$ in the general case.  The
equality (\ref{eq: turaev}) confirms a hypothesis formulated in
\cite{Turaev2}.

The paper  consists of three parts and an Appendix.  The first part  is a review of basic, known facts  about Seifert manifolds.  Its inclusion in the present version of the paper is justified only by   my constant worry to get all the signs right.   The  existent literature  can be quite confusing and/or incomplete about the various orientation conventions.

The second part  deals with the Froyshov invariants. The
computational heart of the paper is \S 2.2 while the applications
are collected in \S 2.3. \S 2.4 contains a number of conjectures
concerning the Froyshov invariants suggested by  numerical
experiments. The most   conceptual one loosely states that  if the
rational homology sphere $N$ is the link of an isolated complex
singularity then the ``most complicated''  negative definite which
bounds $N$ is the minimal resolution of the singularity. The
measure of complexity of a negative definite intersection form is
given by the Elkies invariant described in \S 2.1. The third part
is devoted to the proof of (\ref{eq: cw}) and (\ref{eq: turaev}).

\bigskip

\noindent{\bf Acknowledgements} I want to thank Paolo Lisca  for
his interest in these issues. It was his paper \cite{Lisca}  and
his e-mail questions  on the Froyshov invariants which attracted
my attention to lens spaces. I learned about the Seiberg-Witten
invariants of rational homology spheres from  Weimin Chen who is
one of the pioneers in this subject and I want to thank him for
the useful conversations over the years. As always, Nikolai
Saveliev  has generously shared his knowledge with me. In
particular,  he made  me aware of \cite{BL} which provided the
stimulus for the third part of this paper.   I am indebted to
Frank Connolly  for patiently explaining    the Whitehead torsion
to me and in general, for the many helpful mathematical
conversations.   Finally, I want to thank Yuhan Lim for sending me
his preprints \cite{Lim} and \cite{Lim2}.

\bigskip

\noindent{\bf Orientation conventions} Throughout this paper we will use the following orientation conventions.

\noindent $\bullet$ The boundary of an oriented manifold $M$ is given the outer-normal-first orientation i.e.
\[
{\bf or} (M) = {\rm outer\;\; normal}\wedge {\bf or}(\partial M).
\]
$\bullet$ The total space of a fibration $F\hookrightarrow E\ra B$ is given the fiber-first orientation i.e.
\[
{\bf or}(E) ={\bf or}(F) \wedge {\bf or}(B).
\]
If $E$ is (locally) a principal $S^1$ bundle then the fibers are given the orientation induced by the action of $S^1$.

\noindent $\bullet$  If $\sigma$ is a $spin^c$ structure on an oriented Riemann $3$-manifold $(N,g)$ and ${\bS}_\sigma$ is the associated bundle of complex spinors then  the Clifford multiplication  ${\bc} :\Omega^*(N)\ra {\rm End}\,({\bS}_\sigma)$  is chosen such that
\[
{\bf c}(\ast {\bf 1})=-{\bf Id}.
\]

\bigskip

\tableofcontents

\section{A review of Seifert fibrations}
\setcounter{equation}{0}

The goal of this section is to survey existing results concerning Seifert fibration and, in particular,  clarify the many orientation conventions  concerning the Seifert invariants.

\subsection{Classification results}
In this paper, a {\em Seifert manifold (or fibration}) is a
compact, {\em oriented}, smooth  3-manifold $N$ without boundary,
equipped with an infinitesimally free $S^1$  action.   The orbits
of the  $S^1$-action are called {\em fibers}. A fiber  $S^1\cdot
x$ is called {\em regular} if the stabilizer ${\bf St}_x$ of $x$
is trivial. Otherwise, the fiber is called {\em singular}. In this
case  ${\bf St}_x$ is a cyclic group ${\bZ}_\alpha$ and the order
of this stabilizer  is called the {\em multiplicity} of the fiber.
It is customary to identify ${\bf St}_x$ with the cyclic subgroup
\[
C_\alpha =\{ \exp (\frac{2k\pi \ii}{\alpha});\; k=0,1,\cdots ,\alpha -1\}\subset S^1.
\]
For  brevity set $\rho_\alpha := \exp(\frac{2\pi \ii}{\alpha})$.  The {\em base} of the Seifert fibration is the  space of orbits $\Sigma:= N/S^1$. Topologically, it is a compact oriented surface but smoothly, it is a 2-dimensional orbifolds. The orbifold  singularities are all cone-like and  correspond bijectively to the singular fibers.

Equip $N$ with an $S^1$-invariant Riemann metric $h$. Suppose $F\subset N$ is a singular fiber  of multiplicity $\alpha$ containing the point $x$. The bundle $TN\!\mid_F$ splits orthogonally as
\[
TN\!\mid_F= TF \oplus (TF)^\perp.
\]
Both $TF$ and $(TF)^\perp$ are $S^1$-equivariant bundles over $F$.
The stabilizer $C_\alpha$ of $x$ acts {\em effectively} on $(T_xF)^\perp$. Denote this action by
\[
\tau : C_\alpha \ra {\rm Aut}\,( (T_xF)^\perp).
\]
If we identify $(T_xF)^\perp$ as an oriented vector space with ${\bC}$ then $\tau$  is completely described by an integer $0 < q < \alpha$, $g.c.d.(q,\alpha)=1$  by the formula
\[
\tau(\rho_\alpha)z= \rho_\alpha^qz.
\]
We will denote this action by $\tau_{\alpha, q}$ or, when no
confusion is possible, by $\tau_q$. Following \cite{Or} we call
the pair $(\alpha, q)$ the {\em orbit invariant} of the singular
fiber $F$. Now  denote with $\beta$ the integer uniquely
determined by the requirements
\[
0< \beta < 1, \;\;\;\beta q\equiv 1 \;\;({\rm mod}\;\alpha).
\]
The pair $(\alpha, \beta)$ is called the {\em (oriented, normalized,) Seifert invariant} of the singular fiber $F$.

Using the principal $C_\alpha$-bundle
\[
P_\alpha= (S^1\ra S^1),\;\; z\mapsto z^\alpha
\]
and the representation $\tau_q$ we can form the associated  $S^1$-equivariant  line bundle
\[
E_{\alpha, q}:= P_\alpha \times_{\tau_q} {\bC}\ra S^1.
\]
 The $S^1$-action on $E_{\alpha,q}$ is induced from the obvious action on $S^1\times {\bC}$
\[
e^{\ii \theta} \cdot (z_1, z_2) =(\, e^{\ii\theta}z_1, z_2\,),\;\;|z_1|=1, \;\;z_2\in {\bC}
\]
which commutes with the  action of $C_\alpha$
\[
\rho_\alpha(z_1, z) =(\rho_\alpha z_1, \rho_\alpha^{-q}z_2).
\]
To describe this more explicitly note first that $E_{\alpha,q}$ is
diffeomorphic to $S^1\times {\bC}$. Such an diffeomorphism can be
obtained using the $C_\alpha$ invariant map
\[
T:S^1\times {\bC} \ra S^1\times {\bC}, (z_1, z_2)\stackrel{T}{\mapsto}(\zeta_1, \zeta_2)=(z_1^\alpha, z_1^qz_2).
\]
Then  we can regard $(\zeta_1, \zeta_2)$ as global coordinates on $E_{\alpha, q}$ and we can describe the $S^1$-action  by
\[
e^{\ii\theta}(\zeta_1,\zeta_2)=Te^{\ii\theta}\cdot(z_1,z_2)=(e^{\ii\alpha\theta}\zeta_1, e^{\ii q\theta}z_2).
\]
We have the following result (see \cite{Or}).

\bigskip

\noindent{\bf The Slice Theorem}   {\em There exists an $S^1$-invariant open neighborhood $U$ of $F$ in $N$,   an $S^1$-invariant open neighborhood $V$ of the zero section of $E_{\alpha, q}$ and an $S^1$-equivariant diffeomorphism
\[
\phi: V\ra U
\]
which maps the zero section to $F$ and ${\bf 1}\in S^1$ to a given fixed point $x\in F$.}

\bigskip

Denote $D_r$ denotes the disk of radius $r$ in the fiber of $E_{\alpha, q}$ over ${\bf 1}\in S^1$  i.e.
\[
D_r =\{(1, \zeta_2) \in E_{\alpha,q};\; |\zeta_2|\leq r\}.
\]
The surface $\phi(D_r)$ will be called a {\em slice} of the
$S^1$-action. For simplicity, we will continue to denote it by
$D_r$. Its boundary, equipped with the induced orientation, will
be denoted by $\vec{\sigma}$.  It can be explicitly described by
the parameterization
\[
(\zeta_1,\zeta_2)=(1, re^{\ii t}),\;\;t\in[0,2\pi].
\]

Denote by $\Delta(r)=\Delta_{\alpha, \beta}$ the bundle of disks of radius $r$ determined by $E_{\alpha, q}$ and  set $S(r)=S_{\alpha,\beta}:=\partial \Delta_{\alpha, \beta}$. $\Delta(r)$ is topologically a solid torus.  It is usually known as the {\em fibered torus} corresponding to the Seifert invariants $(\alpha, \beta)$. Endow $S(r)$ with the induced orientation.   $S(r)$ is equipped with a  {\em free}  $S^1$-action. Denote by $\vec{f}$ such an orbit,  endowed with the induced orientation. It can be described explicitly by the curve
\[
(\zeta_1,\zeta_2)=(e^{\ii\alpha t}, e^{\ii qt}),\;\;t\in [0,2\pi].
\]
$\vec{f}$ meets $\vec{\sigma}$ geometrically $\alpha$-times.     In fact, with all the above orientation conventions in place, we also have $\vec{\sigma}\cdot \vec{f}=\alpha$, algebraically as well.

A {\em section} of the $S^1$-action on $S(r)$ is a closed, oriented curve $\vec{s}$ such that $\vec{s}\cdot \vec{f}=1$ both algebraically and geometrically. There exist many sections. We want to show  that there exists a  section  satisfying the  homological condition
\be
\vec{\sigma} =\alpha \vec{s} +\beta \vec{f}.
\label{eq: sec}
\ee
Clearly the above condition uniquely determines the homology class of $\vec{s}$ in $S_r$.

To find a section satisfying (\ref{eq: sec}) we first choose a longitude, i.e. a simple, closed, oriented curve $\vec{\lambda}$ such that $\vec{\sigma}\cdot \vec{\lambda}=1$. There is no unique choice, but two choices differ by a multiple of  $\vec{\sigma}$.  Note that the image of such a $\vec{\lambda}$ in $H_1(\Delta(r), {\bZ})$  coincides with the positive generator, or via $\phi$, with the singular fiber $F$. Then
\[
\vec{f}= u\vec{\sigma} + v \vec{\lambda},\;\;u,v\in {\bZ}
\]
and since $\vec{\sigma}\cdot \vec{f}=\alpha$ we deduce $v=\alpha$ i.e.
\[
\vec{f}= u\vec{\sigma} +\alpha\vec{\lambda}.
\]
Since  $\vec{f}$  ``wraps'' along $\vec{\sigma}$ $q$-times, the
coordinate $u$   is  uniquely determined mod $\alpha$, more
precisely $u\equiv q$ mod $\alpha$. Now  choose $\vec{\lambda}$ so
that $u=q$  i.e.
\be
\vec{f}=q\vec{\sigma} +\alpha \vec{\lambda}_0. \label{eq: fiber}
\ee We call $\vec{\lambda}$ the {\em canonical longitude}. The
sought for section $\vec{s}$ has a decomposition
\[
\vec{s}= x\vec{\sigma} +y\vec{\lambda}
\]
subject to the constraint (\ref{eq: sec}) which becomes
\[
\vec{\sigma}=(x\alpha + \beta q)\vec{\sigma} +(\beta\alpha +\alpha y)\vec{\lambda}.
\]
Since $\beta q \equiv 1$ mod $\alpha$,    there exists an unique
pair  $(x_0,y_0)$ so that the above equality is satisfied. More
precisely
\[
x_0 =(1-\beta q)/\alpha,\;\;y_0=-\beta.
\]
Thus the {\em canonical section}, determined by (\ref{eq: sec}) is
\be
\vec{s}  =x_0 \vec{\sigma} -\beta\vec{\lambda}.
\label{eq: sec1}
\ee

We can now use these notions to describe the structure of Seifert
fibrations.  Suppose the Seifert fibration  has $m\geq 1$ singular
fibers $F_{x_1}, \cdots ,F_{x_m}$ with normalized Seifert
invariants $(\alpha_1, \beta_1)$, ... ,$(\alpha_m, \beta_m)$.
Delete small, pairwise disjoint, $S^1$-invariant neighborhoods
$U_1, \cdots, U_m$ of the singular fibers, determined by  the
Slice Theorem. We get a $3$-manifold with boundary
\[
N'=N\setminus \left( \bigcup_{i=1}^m U_i\right)
\]
equipped with a free $S^1$-action.  This is a principal $S^1$-bundle $S^1\hookrightarrow N'\ra \Sigma':=N'/S^1$. The restriction of this bundle to $\partial \Sigma'$ has canonical sections, determined by (\ref{eq: sec}). In other words, it is trivialized along the boundary.  Such a bundle is completely determined topologically  by an integer $b$,  the  relative degree (or Euler number). Here we have to warn the reader that our $b$ differ by a sign from the conventions in \cite{JN} or \cite{NR}.

We can now reconstruct  $N$ from $N'$ and the equivariant  bundles $E_{\alpha_i, q_i}$ by attaching the fibered torus $\Delta_{\alpha_i, \beta_i}$ to the $i$-th boundary component $\partial_i N'$ of $N'$ using the attaching rules (\ref{eq: sec})
\[
\vec{\sigma}_i = \alpha_i \vec{s}_i +\beta_i \vec{f}_i, \;\;i=1, \cdots, m.
\]
We have to be very careful about the orientation conventions. More precisely,   $(\vec{\sigma}_i, \vec{\lambda}_i)$ and $(\vec{s_i}, \vec{f}_i)$ are  bases of $H_1(\partial \Delta_{\alpha_i,\beta_i}(r), {\bZ})$ compatible with the orientation of $\partial \Delta_{\alpha_i, \beta_i}$    regarded as boundary of $\Delta_{\alpha_i,\beta_i}$. Denote by $\vec{\mu}_i$ the $i$-th boundary component of $\Sigma'$ oriented accordingly. We regard it as an oriented curve on $\partial N'$ via the above trivialization of $N'\!\mid_{\partial \Sigma'}$. Then $(\vec{\mu}_i, \vec{f}_i)$  is compatible  with the orientation of $\partial \Delta_{\alpha_i, \beta_i}$  regarded as a component of $\partial N'$.   On the other hand, $\vec{\mu}_i=-\vec{s}_i$ in $H_1(\partial\Delta_{\alpha_i,\beta_i}, {\bZ})$.   The attaching map
\[
\gamma_i : \partial_i N' \ra \partial \Delta_{\alpha_i, \beta_i}
\]
is given by the identifications (\ref{eq: sec1}) and (\ref{eq: fiber})
\be
\vec{\mu}_i=-\vec{s}_i \mapsto - x_i\vec{\sigma}_i +\beta_i \vec{\lambda}_i,\;\;\vec{f}_i\mapsto q_i \vec{\sigma}_i +\alpha_i\vec{\lambda}_i
\label{eq: glue1}
\ee
where $\alpha_ix_i +\beta_iq_i=1$. If we choose angular coordinates $(\theta^1, \theta^2)$ on $\partial_i N'$   and $(\vfi^1, \vfi^2)$ on $\partial \Delta_{\alpha_i, \beta_i}$ such that $\vec{\mu}_i:=(\theta^1=t,\theta =0)$, and $\vec{s}_i=(\vfi^1= t,\vfi^2 = 0)$, $t\in [0,2\pi]$, then the  above gluing map can be given the matrix description
\[
\left[
\begin{array}{c}
\vfi^1 \\
\vfi^2
\end{array}
\right] =   \left[
\begin{array}{rc}
-x_i & q_i \\
\beta_i & \alpha_i
\end{array}
\right] \cdot \left[
\begin{array}{c}
\theta^1 \\
\theta^2
\end{array}
\right]
\]
The  above matrix has determinant $-1$  and inverse
\be
\Gamma_{\alpha_i,\beta_i} :=
\left[
\begin{array}{rc}
-\alpha_i & q_i \\
\beta_i & x_i
\end{array}
\right].
\label{eq: glue2}
\ee
It is customary to regard the above procedure  the opposite way, as attaching $\Delta_{\alpha_i,\beta_i}$ to $\partial_i N'$ via the orientation reversing map
\[
\Gamma_{\alpha_i, \beta_i}: \partial \Delta_{\alpha_i, \beta_i} \ra \partial_i N'.
\]
Now denote by $\ell$ the rational number
\[
\ell =b-\sum\frac{\beta_i}{\alpha_i}.
\]
It is called the rational number of the Seifert fibration.  The {\em normalized Seifert invariant} of $N$ is defined as the collection
\[
(g, b,  (\alpha_1, \beta_1), \cdots, (\alpha_m ,\beta_m))
\]
where $g$ denotes the genus of $\Sigma$.

The above discussion shows that any Seifert manifold  is uniquely determined (up to an $S^1$-equivariant diffeomorphism) by its Seifert invariant. Moreover,  given a collection as above (with obvious restrictions on the pairs $(\alpha_i, \beta_i)$  one can construct a Seifert manifold with precisely this normalized Seifert invariant.  To  see this, we need only to explain how to construct an $S^1$-bundle over a Riemann surface $\Sigma'$  of genus $g$, obtained from a closed surface $\Sigma$ by deleting  $m$ pairwise disjoint disks $D_1,\cdots, D_m$.  This construction proceeds as follows.

First, delete one more disk $D_0$ from $\Sigma'$ which does not  meet $\partial \Sigma'$. Set $\Sigma''=\Sigma'\setminus D_0$ and $N''=\Sigma''\times S^1$.  Denote by $\partial_0\Sigma''$ the new boundary component.  Now attach $D_0\times S^1$ to $\partial_0 N''$ using the orientation reversing map
\[
\Gamma_b: \partial D_0 \times S^1 \ra \partial_0 N''
\]
given by the matrix
\be
\Gamma_b=\left[
\begin{array}{cc}
-1 & 0 \\
-b & 1
\end{array}
\right]. \label{eq: glue4} \ee The $S^1$-bundle obtained in this
manner is trivialized along the boundary of $\Sigma'$ and has
relative degree $b$. (For $m=0$  this construction mimics the
construction of the holomorphic line bundle on $\Sigma$ associated
to the divisor $b P$, where $P$ is the center of $D_0$.) If we set
$\vec{\sigma}_0=\partial D_0$, $\vec{\lambda}_0=\{1\} \times S^1
\subset \partial D_0 \times S^1 $, $\vec{\mu}_0 = \partial_0
\Sigma''\times \{1\} \subset \partial_0 N''$ and $\vec{f}_0
=\{1\}\times S^1 \subset \partial_0 \Sigma'' \times S^1$, then the
above gluing map produces the identifications
\be
\mu_0 \mapsto -\vec{\sigma}_0 -b\vec{\lambda}_0, \;\; \vec{f}_0\mapsto \vec{\lambda}_0.
\label{eq: glue3}
\ee

Often it is useful to work with un-normalized  Seifert invariants. These are collections
\[
{\bf S}=(g, b, m;(\alpha_1, \beta_1), \cdots , (\alpha_m, \beta_m) )
\]
such that $g.c.d(\beta_i, \alpha_i)=1$, $\alpha_i \neq 0$. Two collections  ${\bf S}$ and ${\bf S}'$ are called equivalent if $g=g'$, the collection of $\alpha_i$-s  not equal to $1$ coincides (including multiplicities) with the collection of $\alpha_j'$-s not equal to $1$ and
\[
b-\sum_i\frac{\beta_i}{\alpha_i}=b'-\sum_i\frac{\beta_j'}{\alpha_j'}.
\]
Clearly, by choosing a section other than the canonical one, we
arrive at an un-normalized Seifert invariant. We refer the reader
to \cite{JN} or \cite{Or}  for a proof of the fact that equivalent
un-normalized Seifert invariants lead to $S^1$-diffeomorphic
Seifert manifolds.

Using the normalized Seifert invariants, and the above gluing description of a Seifert manifold, it very easy to determine its  fundamental group via  Van Kampen's theorem. The   fundamental group of $N''$ has generators
\[
\ a_j, b_j,\vec{\mu}_i, \vec{f},\;1\leq g,\;\;0\leq i\leq m
\]
and relations
\[
[a_1,b_1]\cdots [a_g, b_g] \vec{\mu}_0\cdots \vec{\mu}_m =[a_j,\vec{f}]=[b_j,\vec{f}]=[\vec{\mu}_i,\vec{f}]=1.
\]
Attaching the solid torus $D_0\times S^1$ we introduce a new relation given by (\ref{eq: glue3}) namely
\[
\vec{\mu}_0=\vec{f}^{-b}.
\]
Attaching the fibered torus $\Delta_{\alpha_i, \beta_i}$ we introduce an additional  generator, $\vec{\lambda}_i$ and additional   relations, given  by (\ref{eq: glue1}), namely
\[
\vec{\mu}_i =\vec{\lambda}_i^{\beta_i},\;\;\vec{f}=\vec{\lambda}_i^{\alpha_i}.
\]
Recall that $\vec{\mu}_i$ is a section of $N''$ over $\partial_i\Sigma''$ {\em oriented as a boundary component of $\Sigma''$} and $\vec{f}$ denotes the class of a regular fiber. $\vec{\lambda}_i$  can be expressed in terms of $\vec{\mu}_i$ and $\vec{f}$ by $\vec{\lambda}_i=\vec{\mu}_i^{q_i}\vec{f}^{x_i}$ where $\alpha x_i+\beta_i q_i=1$. Thus, attaching the fibered torus $\Delta_{\alpha_i,\beta_i}$ has the overall effect of introducing the relation
\[
\vec{\mu}_i^{\alpha_i}=\vec{f}^{\beta_i}.
\]
Thus the fundamental group of $N$ can be given the presentation

\medskip

\noindent $\bullet$ {\bf Generators} $a_j, b_j$ ($1\leq j\leq g$), $\vec{\mu}_i$ ($1\leq i\leq m$), $\vec{f}$.

\noindent $\bullet$ {\bf Relations}
\[
\vec{f}^{-b}[a_1,b_1]\cdots [a_g, b_g] \vec{\mu}_1\cdots \vec{\mu}_m =[a_j,\vec{f}]=[b_j,\vec{f}]=[\vec{\mu}_i,\vec{f}]=\vec{\mu}_i^{\alpha_i}\vec{f}^{-\beta_i}=1.
\]

In \cite{FS} the   Seifert manifolds were given a different interpretation in terms of $V$-line bundles over $V$-surfaces. This lead to  different Seifert invariants.  We conclude this subsection with a description of the relationship between the Seifert invariants of \cite{FS} (or \cite{N}) and the Seifert invariants used in this paper.

As explained in \cite{FS},     there is an alternative procedure of obtaining all the Seifert manifolds. Start with a $V$-surface $\Sigma$ with $m$ singular  points $x_1, \cdots, x_m$ with isotropies $C_{\alpha_1}, \cdots ,C_{\alpha_m}$.  Pick a  complex line $V$-bundle    $L\ra \Sigma$ such that the isotropies  in the fibers over the singular points are given  by the representations
\[
\tau_{\alpha_i, \omega_i}: C_{\alpha_i}\ra U(1),\;\; \tau_{\alpha_i, \omega_i}(\rho_{\alpha_i})=\rho_{\alpha_i}^{\omega_i}.
\]
Above, $\omega_i$ are integers satisfying the conditions
\be
0< \omega_i < \alpha_i,\;\;g.c.d(\alpha_i, \omega_i)=1.
\label{eq: norm}
\ee
Then the unit circle bundle $N=S(L)$ determined by $L$ is a Seifert manifold. In \cite{N} we defined the Seifert invariants as the collection
\[
(g, \ell, m;(\alpha_1, \omega_1), \cdots ,(\alpha_m, \omega_m) )
\]
where $\ell$ is the  rational degree of $L$.   We will refer to these as the {\em Seifert $V$-invariants}.

We want to show that the  normalized Seifert invariants (as defined in this paper) of $N$ are
\be
\beta_i:=\alpha_i-\omega_i
\label{eq: on1}
\ee
and
\be
b= \ell +\sum_i\frac{\beta_i}{\alpha_i}.
\label{eq: on2}
\ee
To establish  these facts  we have to understand the orbit invariants of the singular fibers of $S(L)$.

A neighborhood of the singular fiber of $S(L)$ sitting over the singular point  $x=x_i$ can be described  as the $C_{\alpha}$-quotient  of the $C_{\alpha}$-equivariant $S^1$-bundle
\[
T_{\alpha, \omega}\ra D={|z|<1}\subset {\bC}
\]
where $C_\alpha$ acts on $T_{\alpha, \omega}$ by
\[
\rho_\alpha\cdot(z_1, z_2)=(\rho_\alpha z_1, \rho_\alpha^\omega z_2),\;\;(|z_1|<1, \;|z_2|=1)
\]
while $S^1$ acts by
\[
e^{\ii\theta}(z_1,z_2)=(z_1, e^{\ii\theta}z_2).
\]
Note that there exists a natural diffeomorphism
\[
T_{\alpha, \omega}/C_\alpha \ra D\times S^1
\]
induced by the $C_\alpha$-invariant map
\[
A:D\times S^1\ra D\times S^1,\;\;(z_1, z_2) \stackrel{A}{\mapsto}(\zeta_1,\zeta_2)= (z_1z_2^{-s}, z_2^\alpha)
\]
where $s\omega \equiv 1$ mod $\alpha$. We see that $T_{\alpha, \omega}/C_\alpha$ admits an $S^1$-action given by
\[
e^{\ii\theta}(\zeta_1, \zeta_2) = Ae^{\ii\theta}(z_1,z_2)= (e^{-\ii s\theta}\zeta_1, e^{\ii\alpha\theta}\zeta_2).
\]
Hence, the orbit invariants $(\alpha, q)$ of this action are $(\alpha, q)=(\alpha, -s)$.  Thus $\omega q\equiv -1$ mod $\alpha$ so that $\omega\equiv -\beta$ mod $\alpha$. The equality (\ref{eq: on1}) now follows immediately from the normalization condition (\ref{eq: norm}).  We leave the equality (\ref{eq: on2}) to the reader.

We want   to clarify one point. Denote by $|L|$ the desingularization of $L$ (described in \cite{N}).  Then
\be
\deg|L|=\deg L -\sum_i\frac{\omega_i}{\alpha_i}= \ell +\sum_i\frac{\beta_i}{\alpha_i}-m = b-m.
\label{eq: on3}
\ee

The description of Seifert fibrations via line $V$-bundles  has its computational advantages. It allows a very convenient description of the cohomology group $H^2(N, {\bZ})$. We  include it here for later use.

Consider a Seifert fibration  $N$ over a $2$-orbifold $\Sigma$
defined as the unit circle bundle determined by a line $V$-bundle
$L_0\ra \Sigma$.  Suppose  the singularities of $\Sigma$ have
isotropies $\alpha_1, \cdots, \alpha_m$ while the isotropies of
$L_0$ over the singular points are described by
$\rho_{\alpha_i}^{\omega_i}$ as explained above.  Denote by ${\rm
Pic}^t(\Sigma)$ the  space (Abelian group more precisely) of
isomorphism  classes of line $V$-bundles over $\Sigma$. Define  a
group morphism
\[
\tau : {\rm Pic}^t(\Sigma)\ra {\bQ} \oplus {\bZ}_{\alpha_1} \oplus \cdots \oplus {\bZ}_{\alpha_m}
\]
by
\[
\tau(L) = (\deg L, \gamma_1 \;{\rm mod}\; \alpha_1, \cdots , \gamma_m\; {\rm mod}\; \alpha_m)
\]
where $\deg L$ is the rational degree of $L$ and $\gamma_i$ describe the isotropies of $L$ over the singular points of $\Sigma$.  Next, define
\[
\delta: {\bQ}\oplus {\bZ}_{\alpha_1} \oplus \cdots\oplus {\bZ}_{\alpha_m}\ra {\bQ}/{\bZ}
\]
by
\[
\delta(c, \gamma_1, \cdots ,\gamma_m) =\left(c-\sum_i \frac{\gamma_i}{\alpha_i}\right)\;{\rm mod}\, {\bZ}.
\]
In \cite{FS} it is shown that the sequence below is exact
\be
0\ra {\rm Pic}^t(\Sigma) \stackrel{\tau}{\ra} {\bQ} \oplus {\bZ}_{\alpha_1} \oplus \cdots {\bZ}_{\alpha_m}\stackrel{\delta}{\ra} {\bQ}/{\bZ}\ra 0.
\label{eq: picard}
\ee
Moreover, there exists an isomorphism of groups
\be
H^2(S(L_0), {\bZ}) \cong {\bZ}^{2g}\oplus {\rm Pic}^t(\Sigma)/{\bZ}[L_0]
\label{eq: coho}
\ee
where $g$ is the genus of $\Sigma$ and ${\bZ}[L_0]$ denotes the cyclic subgroup of ${\rm Pic}^t(\Sigma)$ generated by $L_0$. The subgroup ${\rm Pic}^t(\Sigma)/{\bZ}[L_0]$  of $H^2(S(L_0), {\bZ})$ consists of  the  Chern classes of the line bundles on $S(L_0)$ obtained  by pullback from line $V$-bundles on $\Sigma$.

\subsection{Plumbing description of Seifert fibrations}
The Seifert manifolds  can be represented as boundaries of certain 4-manifolds naturally determined  from the  Seifert invariant. This is achieved via the  {\em plumbing construction} which we proceed to describe in this section.  Our presentation is greatly inspired from \cite{HNK}. We have to warn the reader that our Seifert invariant conventions differ from those in \cite{HNK}. Ours coincide (except for the sign of $b$) with those in \cite{JN} or \cite{NR}.

We begin by introducing the plumbing construction and a simple way of visualizing it.

Consider two disk bundles $E_1\ra \Sigma_1$  and $E_2\ra \Sigma_2$ over the Riemann surfaces $\Sigma_1$ and $\Sigma_2$, of  degrees $\deg E_1=b$ and $\deg E_1=c$.    These are $4$-manifolds with boundaries circle bundles over Riemann surfaces. As bundles, they are determined by principal $S^1$ bundles over $\Sigma_1$ and $\Sigma_2$ of degrees $b$ and respectively $c$. We present below a computationally friendly way of representing these principal $S^1$-bundles (and thus the associated disk bundles as well). We proceed as follows.

\noindent $\bullet$ Pick a small disk $\Delta_1\subset \Sigma_1$ and   form the trivial $S^1$-bundle $\Delta_1\times S^1\ra \Delta_1$. Set $\Sigma_1'=\Sigma_1\setminus \Delta_1$

\noindent $\bullet$  Fix the integer $r$ and then denote by $N_1'(r)$ the $S^1$-bundle over $\Sigma_1'$ equipped with a trivialization over $\partial\Sigma_1'$ and having relative degree $r$.

\noindent $\bullet$  We can obtain  a degree $b$ $S^1$-bundle over $\Sigma_1$ by attaching the solid torus $\Delta_1\times S^1$ to $\partial N_1'(r)$ via the gluing map $\Gamma_{b-r}:\partial \Delta_1\times S^1\ra \partial N_1'(r)$ described in (\ref{eq: glue4}).

\begin{figure}
\centerline{\psfig{figure=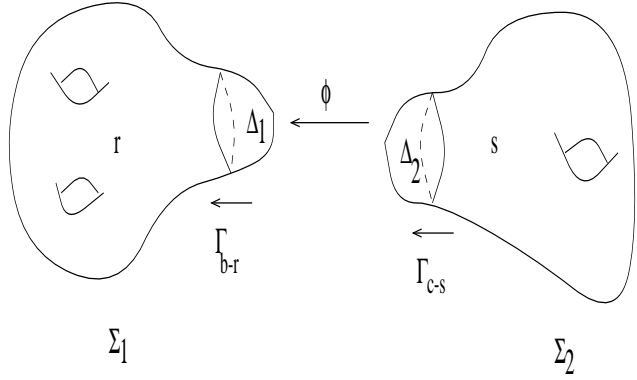,height=2in,width=3.3in}}
\caption{\sl {Plumbing two disk bundles}}
\label{fig: lens1}
\end{figure}

These three steps  can be represented graphically   as in the left-hand side of Figure \ref{fig: lens1}.   We can produce a similar description for  $E_2$ in which $b$ is replaced by $c$ and, instead of $r$, we pick a different integer $s$. This is illustrated in the right-hand-side of Figure \ref{fig: lens1}. We will refer to such a figure  as the {\em diagram of a plumbing}. We will not indicate  the integer $r$ on the diagram when we chose it to be zero. The same convention applies for $s$.

To plumb the disk bundles $E_1$ and $E_2$ proceed as follows.

\noindent $\bullet$  Identify $\Delta_1$ and $\Delta_2$ with the unit disk $D$ in the plane and fix trivializations $E_i\!\mid_{\Delta_i}\ra D\times D$.

\noindent $\bullet$  Now glue $E_2\!\mid_{\Delta_2}$ onto $E_1\!\mid_{\Delta_1}$ using the gluing map
\[
\phi: D\times D\ra D\times D,\;\; (z_1,z_2)\mapsto (z_2,z_1).
\]
The resulting space $E_1\#_\phi E_2$   has apparent corners which can be ``smoothed-out'' to produce a 4-manifold with boundary called the {\em plumbing of the two disk bundles}.      Its boundary  can be alternatively described as follows.

\begin{figure}
\centerline{\psfig{figure=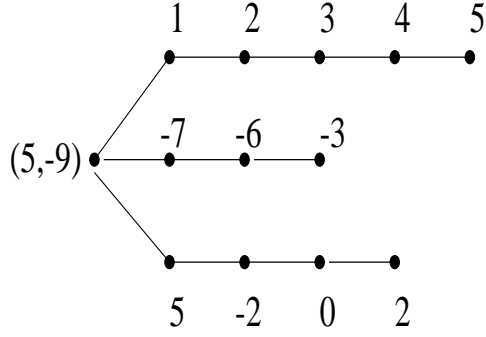,height=1.7in,width=2.5in}}
\caption{\sl {A weighted star-shaped graph}}
\label{fig: lens2}
\end{figure}

\noindent $\bullet$ Attach $\partial N_2'(s)$ to $\partial N_1'(r)$ using the  sequence of gluings
\[
\partial N_2'(s) \stackrel{\Gamma_{c-s}^{-1}}{\ra}\partial \Delta_2\times S^1 \stackrel{\phi}{\ra}\partial \Delta_1\times S^1 \stackrel{\Gamma_{b-r}}{\ra}\partial N_1'(r).
\]
Observe now that $\Gamma_d=\Gamma_d^{-1}$ for any $d\in {\bZ}$. Thus,  the boundary of the plumbing $E_1\#_\phi E_2$ can be obtained by attaching $\partial N_2'(s)$ to $\partial N_1'(r)$  via the gluing map
\[
\Gamma_{b-r}\circ\phi \circ \Gamma_{c-s}.
\]
A {\em star-shaped  graph} is a connected tree with a distinguished vertex  $v_0$ (called the {\em center}) such that the degree of any vertex other than the center is $\leq 2$.  A {\em  branch} of a star-shaped graph is a connected component of a the graph obtained by removing the center.  A {\em weight} on a star-shaped graph $\Gamma$ is a map
\[
w: {\rm Vertex}(\Gamma)\ra ({\bZ}_+\times {\bZ}) \cup {\bZ}
\]
such that
\[
w({\rm center}) \in {\bZ}_+\times {\bZ}
\]
and for any vertex $p\neq {\rm center}$ $w(p)\in {\bZ}$ (see Figure \ref{fig: lens2}).
A weighted star-shaped graph $(\Gamma, w)$  encodes the following topological operations.

\noindent $\bullet$  If the weight of the center is $(g,d)$ associate to it a disk bundle of degree $d$ over  a Riemann surface of genus $g$.

\noindent $\bullet$  To any vertex, other than the center, of weight $n$, associate a degree $n$ disk bundle over $S^2$.

\noindent $\bullet$ Plumb the above disk bundles following the edges of $\Gamma$ i.e. two  bundles are plumbed iff the corresponding edges are  joined by an edge.

In this manner we obtain a $4$-manifold with boundary $P(\Gamma, w)$.  We have the following theorem of von  Randow, \cite{Ran}; see also \cite{Or}.

\begin{theorem}{\rm The boundary of $P(\Gamma, w)$ has a natural structure of Seifert manifold. The Seifert invariants  can be read off the weighted graph $(\Gamma, w)$.}
\end{theorem}

\begin{figure}
\centerline{\psfig{figure=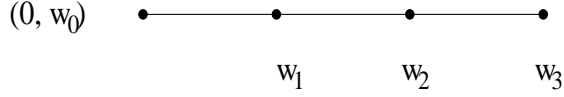,height=.5in,width=2.9in}}
\caption{\sl {A simple plumbing}}
\label{fig: lens3}
\end{figure}

\begin{figure}
\centerline{\psfig{figure=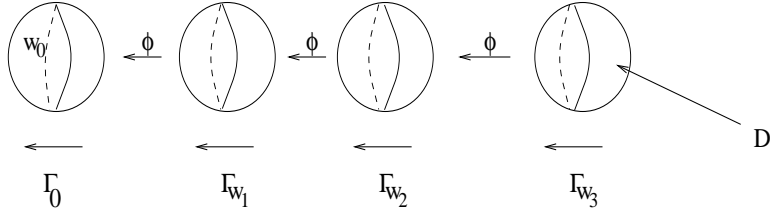,height=1.1in,width=4in}}
\caption{\sl {The diagram of a simple plumbing}}
\label{fig: lens4}
\end{figure}

Let us describe how to read off an {\em un-normalized} Seifert
invariant
\[
(g,b, (\alpha_1, \beta_1), \cdots, (\alpha_m, \beta_m))
\]
of $\partial P(\Gamma, w)$. First of all the number $m$ is
precisely the number of branches of  $\Gamma$. $(g,b)$ is the
weight of the center. Finally, if the weights on the $i$-th branch
are $w_{i1}, \cdots ,w_{ik}$ then the {\em irreducible} fraction
$\alpha_i/\beta_i$ is recovered from the continuous fraction
decomposition
\[
\frac{\alpha_i}{\beta_i}=[w_{i1}, w_{i2}, \cdots ,w_{ik}]
\]
where for any integers $n_1,\cdots ,n_k$ with $n_k\neq 0$ we  define inductively
\[
[n_1, n_2, \cdots, n_k] =n_1 -\frac{1}{[n_2, \cdots ,n_k]} =n_1-\frac{1}{n_2-\frac{1}{n_3 -\frac{1}{\ddots -\frac{1}{n_k}} } }.
\]
We check this on the simple graph depicted in Figure \ref{fig: lens3}.  The boundary of the plumbing is obtained by gluing  the solid disk $D\times S^1$ ($D$ is described in Figure \ref{fig: lens4}) to the boundary of an $S^1$-bundle of  relative degree $w_0$ over a disk. The attaching map  can be read easily from the  diagram in Figure \ref{fig: lens4} and it is
\[
\gamma= \Gamma_0 \circ(\phi \circ \Gamma_{w_1})\circ(\phi\circ \Gamma_{w_2})\circ (\phi\circ \Gamma_{w_3}).
\]
Set
\[
S_b:= \phi \circ \Gamma_b = \left[
\begin{array}{cc}
-b & 1 \\
-1 & 0
\end{array}
\right].
\]

We deduce
\[
\phi\circ \gamma = S_0\circ S_{w_1}\circ S_{w_2}\circ S_{w_3} =\left[
\begin{array}{rc}
-w_2w_3 + 1 &  \ast \\
w_1w_2w_3-w_1-w_3 &  \ast
\end{array}
\right].
\]
Thus
\[
\gamma =\phi\circ(\phi\circ\gamma)=\left[
\begin{array}{rc}
-\alpha & \ast \\
\beta &  \ast
\end{array}
\right]
\]
where $\det \gamma =-1$ and
\[
\alpha =w_1+w_3-w_1w_2w_3,\;\;\beta =1-w_2w_3
\]
Thus $\gamma$  is a gluing map $\Gamma_{\alpha, \beta}$ as in (\ref{eq: glue2})  so that the boundary of this plumbing is the Seifert manifold with un-normalized Seifert invariant
\[
(g=0,b=w_0, (\alpha, \beta)).
\]
Observe now that
\[
\frac{\alpha}{\beta}=\frac{w_1(1-w_2w_3) +w_3}{1-w_2w_3}=w_1 +\frac{1}{\frac{1-w_2w_3}{w_3}}
\]
\[
=w_1 +\frac{1}{-w_2 +\frac{1}{w_3}}=w_1 -\frac{1}{w_2-\frac{1}{w_3}}
\]
as stated in von Randow's theorem.

Observe that there are at least as many plumbing descriptions as
un-normalized Seifert invariants.    In fact there are more
plumbing descriptions than Seifert invariants since the continuous
fraction decomposition  $\frac{\alpha}{\beta}=[w_1,\cdots, w_k]$ is
not unique. E.g. $5/3= [2,1]= [3,1,3]$.    Amongst all continuous
fraction decompositions of a rational number $\alpha/\beta$, ($(\alpha,\beta)=1$, $\alpha>0$) there is a canonical
 one called the Hirzebr\"{u}ch-Jung  plumbing $\lan w_1,\cdots, w_k\ran$ uniquely determined by the requirements
 \[
\lan w_i,\cdots, w_k\ran =w_i-\frac{1}{\lan w_{i+1},\cdots w_k\ran},\;\;{\rm sign}\,(w_i)={\rm sign}\,(\lan w_i,\cdots,
w_k\ran)={\rm sign}\,(\beta)
\]
\[
|w_i|=\lceil\, |\lan w_i,\cdots,w_k\ran| \,\rceil,\;\;\forall i=1,\cdots k
\]
where $\lceil x\rceil$ denotes the smallest integer   $\geq x$.
For example
\[
\frac{8}{5}=<2,3,2>,\;\;
\frac{3}{-2}=\lan-2,-2\ran.
\]
If $N$ is a Seifert manifold  with {\bf normalized} Seifert
invariant $(g,b,m; (\alpha_1,\beta_1),\ldots, (\alpha_m,\beta_m))$
then the {\em Hirzebr\"{u}ch-Jung} plumbing  is obtained from the
Seifert invariant  $(g,b-m,m; (\alpha_i,\beta_i-\alpha_i),\;\;i=1,\ldots,
m)$ using the Hirzebr\"{u}ch-Jung decompositions
\[
\frac{\alpha_i}{\beta_i-\alpha_i}=\lan w_1 ,\cdots w_{k_i}\ran.
\]
Notice that all the weights $w_i$  are negative since
\[
\frac{\alpha_i}{\alpha_i-\beta_i}=\lan -w_1,\cdots , -w_{k_i}\ran.
\]

We conclude with a convention. If the weight of the center of a star-shaped graph  is $(0, w_0)$, that is the associated bundle is a disk bundle over a 2-sphere, then we  say that the plumbing is {\em spherical} and instead of $(0,w_0)$ we will write simply $w_0$.

\subsection{Seifert structures on lens spaces}
We now want to apply the general considerations in the previous subsections to lens spaces.

If $p, q$ are two coprime integers, $p>1$  we define the lens space $L(p,q)$ as the quotient of
\[
S^3:= \{(z_1, z_2)\in {\bC}^2;\; |z_1|^2 +|z_2|^2 =1\}
\]
via the action of $C_p$ given by
\be
\rho_p(z_1,z_2)= (\rho_pz_1, \rho_p^qz_2). \label{eq: action} \ee
Alternatively, we can describe $L(p,q)$ as a result of gluing two
solid tori $D\times S^1$ along their boundaries  using the gluing
map $\Gamma_{q,p}$ (see \cite{JN}).   This shows that we can
regard a lens space  as a Seifert manifold  with (un-normalized)
Seifert  invariant $(g=0, b=0, (q, p))$.   The plumbing discussion
in the previous subsection  shows that we can represent this
Seifert structure  as the boundary of a spherical plumbing given
by a  weighted starshaped graph with one  branch
\[
\stackrel{w_0=0}{\bullet}-\stackrel{w_1}{\bullet}-\cdots -\stackrel{w_k}{\bullet}
\]
where
\[
q/p =[w_1, \cdots ,w_k].
\]
The  {\bf normalized} Seifert invariant  of the above Seifert
fibration is easy to read.  The rational Euler number is $\ell =
-p/q$ so that $b=\lceil-\frac{p}{q}\rceil = -[p/q]$, $\alpha=q$ and
$\beta$ is the remainder of the division $p/q$.  The
Hirzebr\"{u}ch-Jung plumbing corresponding to this  normalized
Seifert in variant is
\[
\stackrel{-w_1}{\bullet}-\cdots -\stackrel{-w_k}{\bullet}
\]
where
\[
\frac{p}{q}=\lan w_1,\cdots, w_k\ran
\]
We will refer to it as the  {\em canonical Hirzebr\"{u}ch-Jung
plumbing} of $L(p,q)$.

In the above graph, each vertex can be regarded  as the center of
another star-shaped spherical graph with possible two branches.
This shows that $L(p,q)$ can be regarded  as a Seifert manifold in
many different ways.  In fact, as explained in \cite{JN} or
\cite{S}, any lens space admits infinitely many Seifert
structures. They all have something in common.  Their bases have
zero genus and they have at most two singular fibers.   Moreover,
as explained in \cite{JN} or \cite{S}, any Seifert fibration over
$S^2$ with at most two singular fibers must be a Seifert fibration
of a lens space.  The  Seifert invariants  of all these Seifert
fibrations are described in  Sec. 4 of \cite{JN}.

  Perhaps, at this point it is instructive to look at the special example of $L(p,1)$. This is the total space of the degree $-p$  circle bundle over $S^2$ and thus has the simple spherical plumbing description
\[
\stackrel{-p}{\bullet}.
\]
This is the canonical Hirzebr\'{u}ch-Jung plumbing of $L(p,1)$. On the other hand
\[
\frac{1}{p}=1-\frac{p-1}{p}=1 -\frac{1}{\frac{p}{p-1}}=1-\frac{1}{1-\frac{1}{1-p}}
\]
so that it has also the plumbing description
\[
\stackrel{0}{\bullet}-\stackrel{1}{\bullet}-\stackrel{1}{\bullet}-\stackrel{1-p}{\bullet}.
\]
If we regard one of the middle vertices as centers we obtain
different   Seifert fibrations structures.    Since $L(p,1)=L(p,
kp+1)$ we can obtain  many other Seifert structures starting from
the continuous fraction decomposition of $(kp+1)/p$.

We will be interested only in those  Seifert structure on a lens
space such that the base is a good orbifold in the sense
described in \cite{S}. This can happen if and only if they have an
(un-normalized) Seifert invariant
\[
(g=0, b=0, (\alpha_1, \beta_1), (\alpha_2, \beta_2) )
\]
satisfying $\alpha_1 =\alpha_2$.  These Seifert structures were determined in \cite{O} for any lens space $L(p,q)$. There are only two of them
\be
{\bf S}_\pm(p,q) =(0,0,(\alpha_\pm , \beta_1^\pm), (\alpha_\pm, \beta_2^\pm))
\label{eq: seiflens}
\ee
which can be explicitly computed as follows.

\noindent $\bullet$ $\alpha_\pm = p/g.c.d.(p, q\pm 1)$

\noindent $\bullet$  $\beta_1^\pm +\beta_2^\pm = \mp g.c.d.(p,
q\pm 1)$.

\noindent $\bullet$
\[
\beta_2^\pm \cdot \frac{q\pm 1}{g.c.d.(p, q\pm 1)} \equiv - 1 \;\;{\rm mod}\;\;\alpha_\pm.
\]
We will refer to the above Seifert structures on $L(p,q)$ as the {\em geometric Seifert structures}. There is a more conceptual description of these structures. To present it,   recall first the Hopf actions of $S^1$ on $S^3$ given by
\[
h_\pm:\;\; (z_1,z_2)\stackrel{e^{\ii\theta}}{\mapsto}(e^{\pm\ii\theta}z_1,e^{\ii\theta}z_2).
\]
The action (\ref{eq: action}) of $C_p$ commutes  with these action of $S^1$ and thus the Hopf actions descend to two infinitesimally free $S^1$-actions on the lens space $L(p,q)$. These define precisely the two geometric Seifert structures.

\subsection{Geometric structures on lens spaces}
All Seifert fibrations  admit   natural geometries, i.e. locally homogeneous Riemann metrics and their universal covers belong to a list of 6 homogeneous spaces (see[S]).    In the case of lens spaces  this geometry is induced from a round metric on their universal cover, $S^3$.    We want to describe those Seifert structures  which interact in a certain way with this metric.

In Sec.1 of \cite{N0} we have described the precise meaning of this interaction (we need a $(K,\lambda)$ structure in the terminology of \cite{N0}). In this case this is equivalent to  asking that the Seifert structures  are the quotient of the Hopf actions on $S^1$ modulo the action (\ref{eq: action}) of $C_p$. In other words, we must restrict to geometric Seifert structures.

Consider a lens space $N=L(p,q)$ equipped with a geometric Seifert structure with (un-normalized) invariant
\[
(g=0,b=0, (\alpha, \beta_1), (\alpha,\beta_2).
\]
The base $\Sigma = N/S^1$ is a $2$-orbifold with at most two conical points of isotropies $C_{\alpha_i}$, $i=1,2$. Denote by $g(R)$ the metric  on $N$ induced by the round metric on the 3-sphere of radius $R$. The radius $R$ will be described below.  The  group $S^1$ acts by isometries of $g(R)$ so that $\zeta$,  the infinitesimal generator of this action, is a Killing vector field.  $\zeta$ is nowhere vanishing and  produces an orthogonal decomposition
\[
TN ={\rm span}(\zeta) \oplus {\rm span}(\zeta)^\perp.
\]
The action of $S^1$ is compatible with this splitting and thus, the metric on ${\rm span}(\zeta)^\perp$ induces an orbifold metric $h$ on $\Sigma$.   Now  fix $R=R_0$ such that
\be
{\rm vol}_h(\Sigma)=\pi. \label{eq: normalize1}
\ee
The radius  $R_0$ can be explicitly determined as follows. Observe
first that the volume of $N$ is equal to
\[
  {\rm length\; regular\; fiber}\times {\rm vol}_h(\Sigma)=
  2\pi^2R_0/p
\]
Since the regular fibers  have length $(1/p)\times${\em(length of a great circle on $S^3(R_0))$} $= 2\pi R_0/p$. Hence
\[
{\rm vol}\,(N)=2\pi^2R_0^2/p.
\]
On the other hand
\[
{\rm vol}\,(N)={\rm vol}\,(S^3(R_0))/p=2\pi^2R_0^3/p
\]
from which we deduce $R_0=1$.

The regular fibers of $N$  are geodesics and have the same length
$2\pi/p$ so  that $\zeta$ has length $1/p$. Denote by
$\vfi\in\Omega^1(N)$ the  $g(R_0)$-dual of  $\zeta$. The metric
$g(R_0)$ can be described as
\[
g(R_0)= \vfi^2 \oplus h.
\]
For $0<r<1$ define the anisotropic rescaling
\[
g_r =(pr)^2\vfi^2 \oplus h.
\]
With respect to this metric the regular fibers have  length $2\pi
r$.  Denote by $\nabla^r$ the  Levi-Civita connection of the
metric $g_1$.    Following \cite{N1} we define  for each $t\in
(0,1]$ an isometry
\[
L_t: (TN,g_{rt}) \ra (TN, g_r), \zeta \mapsto t\zeta,\;\;X\mapsto X\;\;{\rm if}\;\;X\perp \zeta.
\]
Now set
\[
\tilde{\nabla}^{r,t}:= L_t\nabla^{rt}L_t^{-1}.
\]
The connection $\tilde{\nabla}^{r,t}$ is compatible with $g_r$ but
it is not symmetric. In \cite{N1} we have shown that the limit
$\lim_{t\ra 0}\tilde{\nabla}^{r,t}$ exists and defines a
connection compatible with the metric $g_r$. We will call this
limit the {\em adiabatic Levi-Civita} connection  of the metric
$g_r$ and we will denote it by $\tilde{\nabla}^r$.

Observe that a lens space admits two geometric Seifert structures. Arguing as above we obtain two families of Riemann metrics $g_r$ and $h_\rho$.  Both have positive scalar curvature (for $r,\rho \ll 1$) and there exist values $r_0,\rho_0>0$ (which need not be equal) such that the metrics $g_{r_0}$ is homothetic to the metric $h_{\rho_0}$.

\section{Froyshov invariants}
\setcounter{equation}{0}

\subsection{Froyshov's theorem}
For the reader's convenience we include here a brief description
of the Froyshov invariant of a rational homology sphere. For
details we refer to the original source, \cite{Fr}.

Suppose $N$ is rational homology sphere equipped with a Riemann metric $g$. Pick a divergence free $1$-form $\nu$, thought of as a perturbation parameter for the $3$-dimensional Seiberg-Witten equations $SW(g,\nu, \sigma)$ on $(N,g,\sigma)$, where $\sigma$ is a $spin^c$ structure on $N$.  Denote by ${\bS}_\sigma$ the  bundle of complex spinors associated to $\sigma$ and set $\det \sigma =\det {\bS}_\sigma$. The pair  $(g,\nu)$ is said to be {\em good} iff the following hold.

\medskip

\noindent $\bullet$ The  irreducible solutions of $SW(g,\nu,\sigma)$ are nondegenerate for all $\sigma$.

\noindent $\bullet$ If $\theta =(\psi = 0, A_\sigma)$ is the reducible solution of $SW(g,\nu, \sigma)$ then  $\ker \dir_{A_\sigma} =0$ where $\dir_{A_\sigma}$  denotes the Dirac operator  on ${\bS}_\sigma$  coupled with the connection $A_\sigma$ on $\det \sigma$.

\noindent $\bullet$  If nonempty, the  spaces of gradient flow lines (of the  $3$-dimensional Seiberg-Witten energy functional)  which connect irreducible solutions form smooth moduli spaces of the correct dimension.

\medskip

If the pair $(g,0)$ is good then we will simply say the metric $g$ is good.

For any  irreducible solution $\alpha$ of $SW(g,\nu, \sigma)$ denote by $i(\alpha,\theta)$ the  virtual dimension of the space of tunnelings (= connecting gradient flow lines) from $\alpha$ to $\theta$. Define $m=m(g,\nu, \sigma)$ as the smallest nonnegative integer such that there are no tunnelings $\alpha \ra \theta$ with $i(\alpha, \theta)=2m+1$. Now define
\[
{\bf Froy}(N, g,\nu, \sigma) := 8m(g,\nu,\sigma)+4\eta(\dir_{A_\sigma})  +\eta_{sign}(g)
\]
where $\eta_{sign}(g)$ denotes the eta invariant of the odd-signature  operator on $N$ determined by the metric $g$.  In \cite{Fr} it was shown the quantity
\[
{\bf Froy}(N, \sigma) :=\inf \{ {\bf Froy}(N,g,\nu, \sigma);\; (g,\nu)\;\; {\rm is \; good}\}.
\]
Now define the Froyshov  invariant of $N$ by
\[
{\bf Froy}(N):= \max_\sigma {\bf Froy}\,(N,\sigma).
\]
To explain the relevance of this invariant in topology we need to
introduce another, arithmetic  invariant.

Consider  a negative definite integer quadratic form $q$ defined
on a lattice $\Lambda$. Set $\Lambda^\sharp:={\rm
Hom}(\Lambda,{\bZ})$. The quadratic form induces a morphism
\[
I_q:\Lambda\ra \Lambda^\sharp
\]
and since $q$ is nondegenerate  the sublattice $I_q(\Lambda)$ has
finite index $\delta_q$ in $\Lambda^\sharp$.  There exists an
induced   {\em rational} quadratic form $q^\sharp$ on
$\Lambda^\sharp$ by the equality
\[
q^\sharp(\xi_1,\xi_2):=\frac{1}{\delta_q}\lan \, \xi_1\, ,\,
I_q^{-1}(\delta_q\xi_2)\,\ran
\]
where $\lan \,?\, ,\, ?\,\ran: \Lambda^\sharp\times \Lambda \ra
{\bZ}$ denotes the natural pairing.  A vector $\xi\in
\Lambda^\sharp$ is called {\em characteristic} if
\[
\lan \xi,x\ran \equiv q(x,x)\,\mod 2,\;\;\forall x\in \Lambda.
\]
We define the {\em Elkies invariant}  of $q$ by the equality
\[
\Theta(q):= {\rm rank}(q) +\max\{q^\sharp(\xi,\xi);\;\xi\; {\rm
characteristic\;vector\; of}\; q\}
\]
 Note that if $q$ is an even, negative definite form then
\be
\Theta(q)={\rm rank}(q)
\label{eq: even}
\ee
since in this case $\xi=0$ is a characteristic vector.  A result of Elkies (\cite{Elk})  states that if $q$ is a negative definite, {\em unimodular} quadratic form then $\Theta(q) \leq 0$ if and only if $q$ is diagonal.

\begin{theorem}{\bf (Froyshov, \cite{Fr}) } {\rm If $X$ is a smooth, oriented, negative definite $4$-manifold bounding the rational homology sphere $N$ then
\[
\Theta(q_X) \leq {\bf Froy}(N)
\]
where $q_X$ denotes the intersection form of $X$.}
\label{th: froy}
\end{theorem}

\subsection{Computations}

Consider  a lens space $N=L(p,q)$  and fix a geometric Seifert fibration structure on it. The discussion in \S 1.4  shows that the Seifert invariants of this structure  has the form
\[
(g=0, b=0, (\alpha, \beta_1), (\alpha, \beta_2)),\;\;\alpha>0.
\]
More explicitly, this is one of the Seifert structures ${\bf S}_\pm(p,q)$ described  in (\ref{eq: seiflens}).

If we regard $N$ as the unit circle bundle determined by a line $V$-bundle over $\Sigma=S^2(\alpha, \alpha)=N/S^1$ then  we deduce that
\be
\ell:=\deg L_0 =-\frac{\beta_1+\beta_2}{\alpha}
\label{eq: deg}
\ee
and the isotropies of $L_0$ over the singular points are given by
\be
\omega_i  = (-\beta_i) \;\;{\rm mod}\;\; \alpha_i,\;\;\;i=1,2.
\label{eq: iso}
\ee
Above and in the sequel,  for any $x,n\in {\bZ}$  we denote by $x$ mod $n$ the smallest nonnegative integer $\equiv x$ mod $n$.  We want to warn the reader that   when $\alpha =1$ the above Seifert structure has no singular fibers and  $N$ is a genuine smooth $S^1$-bundle over $S^2$  of degree $\ell$.

The canonical line bundle $K_\Sigma$ of $\Sigma$ has  rational degree
\be
\kappa: =-\frac{2}{\alpha}
\label{eq: can}
\ee
so that the rational Euler characteristic is
\be
\chi =-\kappa = \frac{2}{\alpha}.
\label{eq: chi}
\ee

Denote by $\eta_{sign}(r)$ the  eta invariant of the odd signature operator of $N$  equipped with the deformed metric $g_r$ (described in \S 1.4). $\eta_{sign}(r)$ was computed in \cite{O}.  To describe it explicitly we need to introduce the {\em Dedekind-Rademacher sums} defined for the first time by Hans Rademacher in \cite{Ra}. More precisely, for  every pair of  coprime integers $\alpha, \beta$, $\alpha >1$ and any $x, y\in {\bR}$  set
\[
s(\beta,\alpha ; x,y):=\sum_{r=1}^\alpha\left(\left(x +\beta\frac{r+y}{\alpha}\right)\right)\left(\left(\frac{r+y}{\alpha}\right)\right)
\]
where for any $r\in {\bR}$  we set
\[
((r))=\left\{
\begin{array}{rc}
0 & r\in {\bZ} \\
\{q\}-\frac{1}{2} & r\in {\bR}\setminus {\bZ}
\end{array}
\right.  \;\;(\{r\}:=r-[r]).
\]
The sums $s(\beta,\alpha):=s(\beta,\alpha; 0,0)$  are the   Dedekind sums studied in great detail in \cite{HZ} and \cite{RG}.

\be
\eta_{sign}(r)= -{\rm sign}(\ell) + \frac{2\ell}{3}(\chi r^2 -\ell^2 r^4) +\frac{\ell}{3} -4s(\omega_1, \alpha) -4s(\omega_2,\alpha).
\label{eq: sign}
\ee

The canonical  $spin^c$ structure on the orbifold $\Sigma $   (with determinant line bundle $K_\Sigma^{-1}$) determines  by pullback a $spin^c$ structure on $N$ which we denote by $\sigma_0$.  This allows us to bijectively identify the collection of $spin^c$ structures on $L$ with the space of isomorphism classes of complex line bundles. Since $H^2(N ,{\bZ})= {\bZ}_p$ is pure torsion, the discussion at the end of \S 1.1 shows that all the line bundles on $N$ are pullbacks of line $V$-bundles. Thus
\be
{\rm Spin}^c(N)\cong {\rm Pic}^t(\Sigma)/{\bZ}[L_0]
\label{eq: spinc}
\ee
 where  ${\rm Spin}^c(N)$ denotes the space of $spin^c$ structures on  $N$.  If $L$ is a line bundle on $N$ then the $spin^c$ structure $\sigma_0\otimes L$ which corresponds to $L$  has determinant line bundle
\[
\det (\sigma_0\otimes L) = L^{\otimes 2}\otimes \det \sigma_0=L^{\otimes 2}\otimes\pi^*K_\Sigma^{-1}
\]
where $\pi:N\ra \Sigma$ is the natural projection.  The associated bundle of complex spinors is
\[
{\bS}_L= L\oplus L\otimes \pi^*K_\Sigma^{-1}.
\]
In \cite{N} it was shown that, up to gauge equivalence, there is a unique flat connection on $\det \sigma_L$ which we denote by $A_L$.     The Levi-Civita connection of $g_r$ and $A_L$  canonically determine a connection on  ${\bS}_L$ compatible with the Clifford multiplication.   Denote by $\dir_L$ the associated Dirac operator and by $\eta_{dir}(L,r)$ its eta invariant.

The results of \cite{N}  show that for $r$ sufficiently small, the unperturbed Seiberg-Witten equations corresponding to the $spin^c$ structure $L$ have only one  reducible   solution. It is also nondegenerate since the scalar curvature of $g_r$ is positive.  Thus, for  $g_r$ is a good metric (in the sense of Froyshov's theorem) for $r\ll 1$  and since there is no Floer homology  we deduce that
\be
{\bf F}_r(\alpha,\beta_1,\beta_1): =\max \{ 4\eta_{dir}(L, r) +\eta_{sign}(r);\; L\in {\rm Pic}^t(\Sigma)/{\bZ}[L_0]\}
\label{eq: upfroy}
\ee
is an upper bound for the Froyshov invariant ${\bf Froy}(L(p,q))$.

We now show how one can use the results of \cite{N} and \cite{N1} to provide explicit descriptions of
\[
 F_r(L):= 4\eta_{dir}(L, r) + \eta_{sign}(r).
\]
We have to distinguish two cases.

\bigskip

\noindent {\bf A.} $\alpha =1$ so that  $N$ is a degree $\ell$ line bundle over $S^2$ or, as a lens space, $N=L(\ell, -1)=L(|\ell|,|\ell|-{\rm sign}(\ell))$. The signature eta invariant is
\be
\eta_{sign}(r)= -{\rm sign}(\ell) + \frac{2\ell}{3}(\chi r^2 -\ell^2 r^4) +\frac{\ell}{3}
\label{eq: sign1}
\ee
 In this case   there is a unique $spin$ structure on $\Sigma =S^2$ which corresponds to the unique holomorphic square root $K^{1/2}$ of $K_\Sigma$.  This   determines by pullback a   spin structure on $N$  and denote by $\sigma_{spin}$ the  $spin^c$ structure  associated to it.    Then
\[
\sigma_{spin}= \sigma_0\otimes \pi^* K_\Sigma^{1/2}.
\]
For each integer $0\leq k < |\ell|$ we denote by $L_k$ the line bundle of degree $k$ over $\Sigma$ and by $\sigma_k$ the $spin^c$-structure
\[
\sigma_{spin}\otimes \pi^*L_k =\sigma_0\otimes\pi^*(K^{1/2}\otimes L_k).
\]
Also let ${\dir}_k$ denote the Dirac operator on ${\bS}_{\sigma_k}$ determined by the unique flat connection on $\det\sigma_k$ and  denote by $\eta_{dir}(k,r)$ its eta invariant. Then
\[
{\rm Spin}^c(N)=\{ \sigma_k;\;0\leq k < |\ell|\}
\]
and
\[
{\bf F}_r(1,\beta_1,\beta_2) = \max \{ F_r(k):= 4\eta_{dir}(k, r) +\eta_{sign}(r) ;\; 0\leq k< |\ell| \}.
\]
In \cite{N}  we computed the eta invariants, not for the operator ${\dir}_k$, but for the adiabatic Dirac operators $D_k$. These are constructed   using the connection on ${\bS}_{\sigma_k}$ induced  by the adiabatic Levi-Civita connection $TN$ and the flat connection $\det \sigma_k$.  The eta invariant of ${\dir}_k$ can be determined using  variational formul{\ae}  corresponding to  the affine deformation $(1-t){\dir}_k +tD_k$.   The difference
$\eta_{dir}(k,r) -\eta(D_k)$ can be expressed as the sum of a continuous (transgression) term and a discontinuity contribution (spectral flow). The transgression term is  expressed in the second transgression formula of \cite{N1} while the analysis in Sec.4 of \cite{N0} shows  that the spectral  flow  contribution is zero  if $r\ll 1$.   We obtain the following results

\noindent $\bullet$ $k=0$ (use Thm. 2.4 of \cite{N1})
\[
\eta_{dir}(k,r)=\frac{\ell}{6} - \frac{\ell}{6}(\chi r^2 -\ell^2 r^4).
\]
\noindent $\bullet$ $0<k<|\ell|$ (use  the equality (2.22)  and the second transgression formula of \cite{N1})
\[
\eta_{dir}(k,r)= \frac{\ell}{6} +\frac{k^2}{\ell} -{\rm sign}(\ell)k
- \frac{\ell}{6}(\chi r^2 -\ell^2 r^4).
\]
Using (\ref{eq: sign1}) we deduce
\[
F_r(k) = \frac{4}{\ell}k^2 -4{\rm sign}(\ell)k +\ell -{\rm sign}(\ell)
\]
We see that $F_r(k)$ is {\em independent of $r$}!!! Thus  we
\be
{\bf Froy}\,(L(\ell,-1))\leq \max \{ F_r(k);\; 0\leq k <|\ell|\}.
\label{eq: max0}
\ee
We have two subcases

\medskip

\noindent ${\bf A}_1$  $\ell <0$. The maximum  above is
\[
F_r([-\ell/2]).
\]
Thus, when $\ell =-2m$ then
\be
{\bf Froy}(L(2m,1))\leq {\bf F}_r(1,\beta_1,\beta_2) =1,\;\; (2m=\beta_1+\beta_2>0)
\label{eq: froy1}
\ee
and when $\ell =-(2m+1)$ then
\be
{\bf Froy}(L(2m+1,1) )\leq {\bf F}_r(1, \beta_1,\beta_2) = 1 -\frac{1}{2m+1},\;\;(2m+1 =\beta_1+\beta_2 >0).
\label{eq: froy2}
\ee

\medskip

\noindent ${\bf A}_2$ $\ell >0$. In this case the maximum is $F_r(0)=F_r(\ell)=\ell$ so that
\be
{\bf Froy}_r(L(\ell, -1))\leq {\bf F}_r(1,\beta_1, \beta_2)= \ell -1,\;\;(\ell=-(\beta_1+\beta_2)>0).
\label{eq: froy3}
\ee

\bigskip

\noindent ${\bf B}$ $\alpha >1$. The computations are similar in spirit to the ones in case {\bf A} but obviously they are more complex due to the presence of singular fibers.

Let $L\ra N$ be a line bundle over $N=S(L_0)$ and set $\sigma=\sigma_0\otimes L\in {\rm Spin}^c(N)$. To compute $\eta_{dir}(\sigma, r):=\eta_{dir}(L,r)$ we need to determine  the {\em canonical representative} of $L$. This is the  unique line $V$-bundle $\hat{L}=\hat{L}_\sigma\ra \Sigma$ satisfying the conditions
\be
\pi^*\hat{L}\cong L
\label{eq: can1}
\ee
\be
\frac{\kappa -2\deg\hat{L}}{2\ell} \in [0,1).
\label{eq: can2}
\ee
Denote by $\rho=\rho(\sigma)\in [0,1)$ the rational number  sitting in the left-hand-side of (\ref{eq: can2}) and by $0\leq \gamma_i=\gamma_i(\sigma) <\alpha$, $i=1,2$ the isotropy of the fibers of $\hat{L}_\sigma$ over the singular points.  Finally set
\[
d(\sigma)= \frac{\kappa}{2} -\ell \rho(\sigma) =\deg\hat{L}_\sigma.
\]
In Proposition 1.10 of \cite{N} we computed the eta invariant  for the adiabatic  Dirac operator $D_L =D_\sigma$   defined  by using the adiabatic connection on ${\bS}_\sigma$ and the flat connection on $\det \sigma$.  To recover the eta invariant of $\dir_\sigma:=\dir_L$ we use a deformation argument as in Case {\bf A} and  we deduce the following results.

\medskip

\noindent $\bullet$ If $\rho(\sigma)=0$ then
\be
\eta_{dir}(\sigma,r) = \frac{\ell}{6} -2 \sum_{i=1}^2 s(\omega_i, \alpha; \gamma_i(\sigma)/\alpha,0) -\sum_{i=1}^2\LL\frac{q_i\gamma_i(\sigma)}{\alpha}\RR - \frac{\ell}{6}(\chi r^2 -\ell^2 r^4)
\label{eq: dirfl0}
\ee
where $0\leq q_i <\alpha$ denotes the inverse of $\omega_i$ mod $\alpha$.

\noindent $\bullet$ If $\rho(\sigma)>0$ then
\[
\eta_{dir}(\sigma, r)=(1-\frac{1}{\alpha}) (1-2\rho) -\ell\rho(1-\rho) +2\rho +\frac{\ell}{6}
\]
\be
-2\sum_{i=1}^2 s(\omega_i, \alpha; \frac{\gamma_i(\sigma) +\omega_i\rho}{\alpha}, -\rho) -\sum_{i=1}^2\left\{\frac{q_i\gamma_i(\sigma) +\rho}{\alpha}\right\}- \frac{\ell}{6}(\chi r^2 -\ell^2 r^4)
\label{eq: dirfl1}
\ee
where $\{x\}$ denotes the fractional part of the real number $x$.

The above formul{\ae} may seem hopelessly useless.      Fortunately, the Dedekind-Rademacher  sums satisfy a reciprocity law (see \cite{Ra})  which makes them  computationally very friendly.  We include here the reciprocity law for later use in this paper. To formulate it we must distinguish two cases.

\noindent  $\bullet$  Both $x$ and $y$ are integers. Then
\be
 s(\beta, \alpha; x,y) +s(\alpha, \beta; y,x) =-\frac{1}{4}+\frac{\alpha^2+\beta^2+1}{12\alpha \beta}.
\label{eq: rec1}
\ee
$\bullet$  $x$ and/or $y$ is not an integer. Then
\be
s(\beta,\alpha; x, y)+s(\alpha, \beta; y,x)=((x))\cdot ((y)) + \frac{\beta^2\psi_2(y) +\psi_2(\beta y+\alpha x) +\alpha^2\psi_2(x)}{2\alpha \beta}
\label{eq: rec2}
\ee
where $\psi_2(x):= B_2(\{x\})$ and $B_2(z)$ is the second Bernoulli polynomial
\[
B_2(z)= z^2-z+\frac{1}{6}.
\]
 Denote by $R(\beta, \alpha; x,y)$ the right hand side  in the above reciprocity identities. Note that $R(\alpha,\beta; y,x)=R(\beta, \alpha; x,y)$.

The reciprocity law, coupled with the identities
\be
s(\beta, \alpha;x,y)=s(\beta -m \alpha, \alpha ; x +my,y),\;\;\forall m\in {\bZ}
\label{eq: rec3}
\ee
reduces the computation of any Dedekind-Rademacher sum to the special case $s(\beta,1;x,y)$ which is
\be
 s((\beta, 1;x,y)=(\beta y +x))\cdot ((y))
\label{eq: rec4}
\ee

Using (\ref{eq: sign}), (\ref{eq: dirfl0}) and (\ref{eq: dirfl1}) we conclude  that when $\rho(\sigma)=0$ we have
\[
F_r(\sigma)=\ell -{\rm sign}(\ell)-8 \sum_{i=1}^2 s(\omega_i, \alpha; \gamma_i(\sigma)/\alpha,0)
\]
\be
-4\sum_{i=1}^2\LL\frac{q_i\gamma_i(\sigma)}{\alpha}\RR-4\sum_{i=1}^2s(\omega_i,\alpha)
\label{eq: froy5}
\ee
and when $\rho(\sigma)>0$ we have
\[
F_r(\sigma)=\ell -{\rm sign}(\ell)+ 4(1-\frac{1}{\alpha}) (1-2\rho) -4\ell\rho(1-\rho) +8\rho
\]
\be
 -8\sum_{i=1}^2 s(\omega_i, \alpha; \frac{\gamma_i(\sigma) +\omega_i\rho}{\alpha}, -\rho) -4\sum_{i=1}^2\left\{\frac{q_i\gamma_i(\sigma) +\rho}{\alpha}\right\}-4\sum_{i=1}^2 s(\omega_i, \alpha)
\label{eq: froy4}
\ee
Note  again the $r$ has disappeared!!!

\begin{remark}{\rm Define more generally, for any metric  $g$ on $L(p,q)$
\[
F=F_{g}: {\rm Spin}^c(L(p,q))\cong H^2(L(p,q),{\bZ})\cong {\bZ}_p\ra {\bQ},\;\;\sigma\mapsto 4\eta_{dir}(\sigma, g)+\eta_{sign}(g).
\]
Note that  $F$ is unchanged  by rescaling the metric
\[
F_{\lambda^2g}=F_{g}
\]
because the eta invariants are invariant under such changes.  Moreover, for the metrics $g_r$ associated to a geometric Seifert structure we have shown that $F_{g_r}$ is independent of $r$.     There are two geometric  Seifert structures ${\bf S}_\pm$ on $L(p,q)$ and correspondingly, two families  of metrics $g_r^\pm$ and thus two functions
\[
F^\pm=F_{g_r^\pm}: {\bZ}_p\ra {\bQ}.
\]
As explained in \S 1.4    there are  two radii $r_\pm$ such that
the metrics $g^\pm_{r_\pm}$ are homothetic. Thus, the two
functions  $F^\pm$ must be equal.  This corresponds to  a
collection of $p$ identities between Dedekind-Rademacher sums.
Numerical  experimentations agree beautifully with this simple
observation.}
\end{remark}

To  put the formul{\ae} to work we need to have a complete list of
the canonical  representatives of the line bundles on $N$.   Given
the isomorphism (\ref{eq: picard}) this reduces to an elementary
number theoretic problem.

 According to (\ref{eq: picard}) any line $V$-bundle on $\Sigma$ can be uniquely represented as a collection
\[
(\frac{i}{\alpha}, j \; {\rm mod}\; \alpha, (i-j)\;{\rm mod} \;\alpha), \;\;i,j\in {\bZ}.
\]
Set $n=(\beta_1+\beta_2)$ so that $\ell =-n/\alpha$. A collection
as above is the canonical representative of a line bundle as above
if
\[
\frac{\kappa -2i/\alpha}{-2n/\alpha}=\frac{i+1}{n} \in [0,1).
\]
Thus, when ${\rm sign}(n) =-1$   we deduce that the complete list
of canonical representatives is
\be
{\cal  R}_n = \{ (\frac{i}{\alpha}, j \; {\rm mod}\; \alpha,
(i-j)\;{\rm mod} \;\alpha);\; i=-1,-2,\cdots, -|n|, \;\;0\leq j
<\alpha\} \label{eq: rep1} \ee while when ${\rm
sign}(n)=1$ the complete set of canonical representatives is
\be
{\cal  R}_n = \{ (\frac{i}{\alpha}, j \; {\rm mod}\; \alpha, (i-j)\;{\rm mod} \;\alpha);\; i=-1,0,\cdots, |n|-2, \;\;0\leq j <\alpha\}.
\label{eq: rep2}
\ee
The invariant $\rho$ of a canonical representative $\nu=(i/\alpha, j, i-j)\in {\cal R}$ is
\be
\rho(\nu) =\frac{i+1}{n}.
\label{eq: ro}
\ee
Notice that we can identify
\[
 I_{n,\alpha}:=\{-1,0, \cdots, |n|-2\}\times {\bZ}_\alpha \sim {\cal
 R}_n
\]
via the correspondence
\[
 (k, j\,\mod \alpha) \sim \nu\mapsto (\frac{{\rm sign}(n)k-c}{\alpha}, j, -{\rm sign}(n)k-c-j)
\]
where $c:=1-{\rm sign}\,(n)$. The functions $\rho,\gamma_1,\gamma_2:{\cal R}\ra {\bQ}$ can now be regarded  as functions on $I_{n,\alpha}$. More precisely
\be
\rho(k,j\, \mod \alpha)=\frac{k+1}{|n|}
\label{eq: ro1}
\ee
and
\be
\gamma_1(k,j\,{\rm mod}\,\alpha)=j,\;\;\gamma_2(k,j\,{\rm mod}\, \alpha)={\rm sign}\,(n)k-c-j.
\label{eq: gamma}
\ee
Finally we can now regard $F_r$ as a function
\[
F_r=F_r(k,j): I_{n,\alpha}\ra {\bQ}
\]
given by (\ref{eq: froy5}), (\ref{eq: froy4}), (\ref{eq: ro1}) and (\ref{eq: gamma}).
Hence
\be
{\bf Froy}(L(p,q))\leq {\bf F}_r(\alpha,\beta_1,\beta_2)= \max_{(k,j)\in I_{n,\alpha}}F_r(k,j).
\label{eq: froy6}
\ee

From the proof of Proposition 7 in \cite{Fr} we deduce that, since  the metrics $g_r$ have positive scalar curvature, our upper estimates are optimal.  We have thus established the following result.

\begin{theorem}{\rm If the lens space $L(p,q)$ is given a geometric Seifert structure with Seifert invariant
\[
(g=0,b=0, (\alpha, \beta_1), (\alpha, \beta_2))
\]
then
\[
{\bf Froy}(L(p,q))= \max_{(k,j)\in I_{n,\alpha}}F_r(k,j).
\]
where the quantities  $F_r(k,j)$ are described by (\ref{eq: froy5})-(\ref{eq: gamma}) if $\alpha >1$ and by (\ref{eq: froy1})-(\ref{eq: froy3}) if $\alpha =1$.}
\label{th: liviu}
\end{theorem}

From (\ref{eq: froy1})-(\ref{eq: froy3}) we deduce that (for $p>0$)
\be
{\bf Froy}(L(p,1))=\left\{
\begin{array}{rlc}
1&,&p=2k \\
1-\frac{1}{p}&,& p=2k+1
\end{array}
\right.
\label{eq: liviu0}
\ee
and
\be
{\bf Froy}(L(p,p-1))=p-1.
\label{eq: liviu1}
\ee

The above theorem reduces the computation  of the Froyshov invariant to an elementary, albeit complex, arithmetic problem.

\begin{ex}{\rm We want to illlustrate the strength  and limitations of this theorem by computing the Froyshov invariants  of $L(p,2)$ $p$ odd and $>3$.   The case $L(3,2)=L(3,-1)$ corresponds to a degree $S^1$-bundle over $S^2$ and we have already dealt with it.

 We will use the invariants ${\bf S}_-(p,2)$ which seem to be computationally friendlier.  We have
\[
\alpha =p,\;\;\ell =-\frac{1}{p},\;\;n=-1,
\]
\[
\beta_2 = p-1,\;\;\beta_1 = 2-p,\;\; \omega_1= 1, \;\;\omega_2=p-2
\]
\[
q_1=1,\;\;q_2 =\frac{p-1}{2}.
\]
Thus
\[
I_{n,\alpha}= \{0\}\times {\bZ}_p
\]
and $\rho: I_{n,\alpha}\ra {\bQ}$ is identically zero.  Moreover, since $k=0$ for all $(k,j)\in I_{n,\alpha}$ we have
\[
\gamma_1(k,j)=j,\;\;\gamma_2 (k,j)=p-j-1,\;\;0\leq j<p.
\]
Hence
\[
F_r(k,j)= F_r(0,j)=-\frac{1}{p} +1 -8s(1,p,;j/p,0)-8s(p-2, p; (p-j-1)/p,0)
\]
\[
-4((j/p))-4 \LL\,\frac{(p-1)(p-j-1)}{2p}\,\RR -4s(1,p)-4s(p-2,p).
\]
We distinguish three cases.

\medskip

\noindent $\bullet$ $0<j<p-1$ The equality can be slightly simplified using the elementary identities
\[
s(p-2, p, (p-j-1)/p,0)= s(p-2,p; -(j+1)/p,0)= -s(2,p; (j+1)/p, 0)
\]
\[
s(p-2,p)=-s(2,p)
\]
\[
-4\LL\,\frac{j}{p}\,\RR=2-\frac{4j}{p}
\]
and
\[
4\LL\,\frac{(p-1)(p-j-1)}{2p}\,\RR = -4\LL\, \frac{p-1}{2}\frac{j+1}{p}\,\RR=-4\LL\,\frac{j+1}{2}-\frac{j+1}{2p}\,\RR
\]
\[
=\epsilon -\frac{2(j+1)}{p}
\]
where $\epsilon =2$ if $j$ is odd and $=0$ if $j$ is even. We deduce
\[
F_r(j) =1-\frac{1}{p} +8(s(2,p;\frac{j+1}{p},0)-s(1,p;\frac{j}{p},0))
\]
\[
+4(s(2,p)-s(1,p)) -\frac{4j}{p}-\frac{2(j+1)}{p}+2+\epsilon.
\]
To compute the Dedekind-Rademacher sums we use the reciprocity formul{\ae} (\ref{eq: rec1})-(\ref{eq: rec4}). We deduce
\[
s(2,p, (j+1)/p,0) =R(2, p; (j+1)/p,0) -s(p,2;0, (j+1)/p)= R(2, p; (j+1)/p,0)
\]
\[
=\frac{5}{4p}\psi_2(0)+\frac{p}{4}\psi_2(\frac{j+1}{p})= \frac{5}{24 p}+\frac{p}{4}\left(\frac{(j+1)^2}{p^2}-\frac{j+1}{p}+\frac{1}{6}\right).
\]
We deduce similarly
\[
s(1,p;j/p,0) = \frac{1}{6p}+\frac{p}{2}\left(\frac{j^2}{p^2}-\frac{j}{p}+\frac{1}{6}\right).
\]
After some elementary manipulations we deduce
\[
8(s(2,p,(j+1)/p,0)-s(1,p;j/p,0))= -\frac{2j^2}{p}+2j\left(1+\frac{2}{p}\right) -2+\frac{7}{3p}-\frac{p}{3}.
\]
A similar argument  leads to the equality
\be
4(s(2,p)-s(1,p))=-\frac{p}{6}+\frac{1}{6p}.
\label{eq: ded2}
\ee
Together, all of the above  yield  after some elementary but tedious computations
\[
F_r(j)=-\frac{2j^2}{p}+2j\left(1-\frac{1}{p}\right)-\frac{1}{2p}-\frac{p}{2} +1 +\epsilon.
\]
The above expression is quadratic in $j$. Its maximum on the {\em discrete} interval  $(0,p-1)$ is achieved for $j$ equal to one of the odd integers closest to the midpoint $(p-1)/2$.   When $p=4k+3$  there is only one such integer $(p-1)/2$ and we deduce
\[
F_r(\frac{p-1}{2})= 2
\]
while when $p=4k+1$ then $j=(p+1)/2$ is a maximum point
\[
F_r(\frac{p+1}{2})= 2-\frac{2}{p}.
\]

\medskip

\noindent $\bullet$ $j=0$ so that
\[
F_r(0)=1-\frac{1}{p} + 8(s(2,p;1/p,0)-s(1,p)) +4(s(2,p)-s(1,p)) +4\LL\frac{1}{2}-\frac{2}{p}\RR
\]
\[
\stackrel{(\ref{eq: ded2})}{=}1-\frac{5}{6p}-\frac{p}{6}-\frac{8}{p} +8(s(2,p;1/p,0)-s(1,p)).
\]
The above Dedekind-Rademacher sums can be computed using the  reciprocity law and, as before, we deduce
\[
8(s(2,p;1/p,0)-s(1,p)) +4(s(2,p)-s(1,p))=\frac{7}{3p}-\frac{p}{3}.
\]
It is now clear that  $F_r$ cannot have a global maximum at $j=0$.  The case $j=p-1$ can be disposed of similarly and we leave it to the reader.

We have shown}
\be
{\bf Froy}(L(p,2))=\left\{
\begin{array}{rlc}
2&,& p=4k-1\\
2-\frac{2}{p}&,&p=4k+1
\end{array}
\right. .
\label{eq: liviu2}
\ee
\label{ex: liviu}
\end{ex}

The above example  suggests that for large $p,q$ the  computational complexity can be overwhelming. On the other hand, these computations can be performed easily with any computer algebra system and, because of the reciprocity law, one can manipulate quite large numbers. Here are the results of some {\it MAPLE} experiments.
\be
{\bf Froy}(L(p,3))=
\left\{
\begin{array}{rlc}
 3&,& p=6k-2\\
2-\frac{2k}{6k-1}&,&p=6k-1\\
3-\frac{3}{6k+1}&,&p=6k+1 \\
2-\frac{k+1}{3k+1}&,&p=6k+2
\end{array}
\right.
\label{eq: liviu3}
\ee
\be
{\bf Froy}(L(p,4))=\left\{
\begin{array}{rlc}
4 &,& p=8k-3 \\
2&,& p=8k-1 \\
4-\frac{4}{p}&,&p=8k+1\\
2-\frac{4}{p}&,& p=8k+3
\end{array}
\right.
\label{eq: liviu7}
\ee
\be
{\bf Froy}(L(p,5))=\left\{
\begin{array}{rlc}
5 &, &p=10k-4 \\
3-\frac{4k-1}{10k-3}&,& p=10k-3 \\
3-\frac{2k}{5k-1}&,&p=10k-2\\
2-\frac{2k}{10k-1}&, &p=10k-1\\
5-\frac{5}{p}&,&p=10k+1\\
3-\frac{2k+2}{5k+1}&,&p=10k+2\\
3-\frac{4k+3}{10k+3}&,&p=10k+3\\
2-\frac{k+2}{5k+2}&,&p=10k+4
\end{array}
\right.
\label{eq: liviu8}
\ee
In \cite{CH} it was shown that the  lens spaces $L(p^2,p+1)$  bound  rational homology balls. Their  Froyshov invariants are
\be
{\bf Froy}(L(p^2, p+1))=\left\{
\begin{array}{rlc}
p+1 &,& p\;{\rm even} \\
p+1-\frac{p+1}{p^2}&,& p\;{\rm odd}
\end{array}
\right.
\label{eq: liviu9}
\ee

\subsection{Topological applications}
    Let us introduce some terminology.  By a {\em special} manifold we will understand a smooth, oriented, even, negative definite $4$-manifold (with or without boundary).  A {\em very special} manifold is a special manifold $X$ such that $H_1(X, {\bZ})$ has no $2$-torsion. The following is a consequence of Froyshov's theorem  \ref{th: froy} coupled with the equality (\ref{eq: even}).

\begin{corollary}{\rm If $N$ is rational homology sphere and $X$ is a special  $4$-manifold bounding $N$ then
\[
b_2(X)\leq {\bf Froy}(N).
\]
In particular, if ${\bf Froy}(N) <1$ then there are no special manifolds bounding $N$.}
\label{cor: 1}
\end{corollary}

Using (\ref{eq: liviu0}), (\ref{eq: liviu2})-(\ref{eq: liviu8}) we deduce immediately the following topological consequence.

\begin{corollary}{\rm (a) The lens  spaces $L(2k+1, 1)$ bound no special manifold.

\noindent (b) If $X$ is a special manifold which bounds one of the spaces  $L(2k,1)$, $L(4k+1,2)$, $L(6k-1,3)$, $L(6k+2,3)$, $L(8k+3,4)$, $L(10k-1,5)$, $L(10k+4,5)$ must have $b_2(X)= 1$. In particular, the intersection form of $X$  is diagonal.}
\label{cor: 2}
\end{corollary}

Part  (a) of this corollary  is surprising because the lens spaces $L(2k+1,1)$ do bound smooth, even $4$-manifolds. Also, observe that $L(2k,1)$ is the total space of the degree $-2k$ circle bundle   over $S^2$ and bound a special manifold, the associated disk bundle $D_{-2p}$ so, in this  case, part (b) of the corollary is optimal.

Notice that
\[
\frac{4k+1}{2}= 2k-\frac{1}{-2}=[2k,-2]
\]
\[
\frac{6k-1}{6k-4}=[\underbrace{2,2,\cdots,2}_{2k-1},-2]
\]
\[
\frac{8k+3}{4}=[2k,-2,-2,-2]
\]
and
\[
\frac{10k-1}{10k-6}=[\underbrace{2,2,\cdots,2}_{2k-1},-4].
\]
We can now use   Recipe (7.10) of \cite{HNK}  to compute the
Rohlin invariants of the lens spaces in Corollary \ref{cor: 2}
(b). We have
\[
\mu(L(4k+1,2)) =0\;\;{\rm mod}\; 16 {\bZ}
\]
\[
\mu(L(6k-1,3)) =-\mu(L(6k-1,6k-4))= -(2k-2)\;\;{\rm mod}\;16{\bZ}
\]
\[
\mu(L(8k+3,4))=-2\;\;{\rm mod}\; 16{\bZ}
\]
and
\[
\mu(L(10k-1,5))=-\mu(L(10k-1,10k-6))= (2-2k)\;\;{\rm mod}\; 16 {\bZ}
\]
Using the definition of the Rohlin invariant we deduce that if $X$
is a very special  $4$-manifold bounding  one of the above lens
spaces,  then its signature $(=-b_2)$ is congruent modulo $16$ to
$\mu$.  On the other hand,  we know from  Corollary \ref{cor: 2}
(b)  that  this signature must be $-1$.   We can draw the
following conclusion.

\begin{corollary}{\rm   There exists no very special manifold $X$  which bounds one of the lens spaces in the list below}
\be
L(4k+1,2), \;L(6k-1,3), \; L(8k+3,4),\; L(10k-1, 5).
\label{eq: list}
\ee
\end{corollary}

We leave the reader formulate other corollaries of the same nature. We would like to present another consequence of a slightly different nature.  It relies on a  recent result of Paolo Lisca.

\begin{theorem}{\bf (\cite{Lisca})}{\rm  Let $(X,\omega)$ be a $4$-manifold with contact boundary equipped with a compatible  symplectic form. Suppose that a connected component of the boundary of $X$ admits a metric with positive  scalar curvature. Then,  the boundary of $X$ is connected and  $X$ is negative (semi)definite.}
\label{th: lisca}
\end{theorem}

It follows from the above theorem that any even, symplectic $4$-manifold,   with contact boundary a lens space, must be special. We have the immediate consequence.

\begin{corollary}{\rm The lens space $L(2k+1,1)$ cannot be the contact boundary of any sympectic manifold with even intersection form.  Also, none of the spaces in the list (\ref{eq: list}) can be the contact boundary of a symplectic manifold with no $2$-torsion in $H_1$.}
\label{cor: 3}
\end{corollary}

\subsection{Some conjectures and speculations}
The  examples  discussed so far suggest that the following
arithmetic conjecture is plausible.

\medskip

\noindent{\bf Conjecture 1}  {\em Suppose $p,q$ are two coprime
integers such that $p>q>1$.  Denote by ${\cal R}_q$ the set of
integers $0\leq u \leq 2q$ such that $g.c.d.(u,q)=1$   Then

\noindent (a) ${\bf Froy}(L(p,q))\leq q$.

\noindent (b) ${\bf Froy}(L(2qk+1-q,q))=q-1$.

\noindent (c) For each $u\in {\cal R}_q$ there exist integers
$A_u$, $B_u$ such that}
\[
{\bf Froy}(L(u+2kq,q))=\frac{A_uk+B_u}{u+2qk},\;\;\forall k\in
{\bZ}_+
\]

\medskip

  If true, part(c) of this conjecture would provide a very fast way of computing  the Froyshov invariants of infinite families of  lens spaces. The pair $(A_u,B_u)$ can be viewed as defining an universal function
\[
 {\cal R}_q\ra {\bZ}^2, \;\; u\mapsto (A_u, B_u).
\]

Numerical experimentations  have  displayed and interesting
phenomenon.   First let us introduce an equivalence relation on
the space of negative definite integral quadratic forms.  Two such
forms $q_1$ and $q_2$ are said to be equivalent if there exist two
unimodular, negative definite diagonal forms $\delta_1,\delta_2$
such that
\[
q_1\oplus \delta_1 \cong q_2\oplus \delta_2.
\]
We denote this equivalence relation by $\sim$ and the set of its
equivalence classes by ${\cal Q}$. Also we denote  by ${\cal Q}_1$ the subset
of equivalence classes containing unimodular forms. Since the
Elkies invariant of an unimodular diagonal form is zero we deduce
that $\Theta$ defines a map
\[
\Theta: {\cal Q} \ra {\bQ}
\]
such that $\Theta({\cal Q}_1)\in 8 {\bZ}_+$. It is believed (see
\cite{Elk} and \cite{Gaulter}) that  (at least  for unimodular
forms) $\Theta$ provides a ``measure of complexity''  of a
negative definite quadratic form i.e. for any $k\in {\bZ}_+$ the
set
\[
\{q\in {\cal Q}_1;\; \Theta(q) \leq 8k\}
\]
is finite. (In \cite{Gaulter}   this result is proved for $k<4$).

Thus, Froyshov result can be (loosely) interpreted as  describing
a topological upper bound  for the complexity of the negative
definite manifolds with boundary a given rational homology sphere.

On the other hand, the lens spaces,  are links of complex surface
singularities. More concretely, they are links of   quotient
singularities.    These   singularities can be resolved   and the
effect is a complex, negative definite manifold  with (oriented)
boundary the given   lens space.   There is a canonical way of
performing such a resolution introduced by Hirzebr\"{u}ch  (see \cite{BPV}) and
topologically, this  resolution coincides with the canonical
Hirzebr\"{u}ch-Jung plumbing associated  to the given lens space.
Denote by $HJ(p,q)$ the Hirzebr\"{u}ch-Jung plumbing of  $L(p,q)$,
by $S_{p,q}$ the intersection form of $HJ(p,q)$ and
 by $\Theta_{p,q}$ the  Elkies invariant  of  $S_{p,q}$.   Froyshov's theorem  guarantees  that $\Theta_{p,q}\leq
{\bf Froy}(L(p,q))$. We claim the following  stronger result is
true.

\bigskip

\noindent {\bf Conjecture 2} $\Theta_{p,q}= {\bf Froy}(L(p,q))$.

\bigskip

Here is some evidence supporting this conjecture.

\begin{proposition} {\bf Conjecture 2} {\rm  is true for
$L(p,1)$ and $L(p,2)$}
\label{prop: conj2}
\end{proposition}

\noindent {\bf Proof} (a) $L(p,1)$.   We know that
\[
HJ(p,1)\cong  \stackrel{-p}{\bullet}
\]
with intersection form $S_{p,1}=(-p)$.  When $p$ is even, this intersection
form is even so that
\[
\Theta_{p,1}= 1\stackrel{(\ref{eq: liviu0})}{=} {\bf
Froy}(L(p,1)).
\]
When $p$ is odd  then
\[
\Theta_{p,q}=1+\max\{ \lan S^{-1}_{p,1}x,x\ran\, ;\;\;x\in
{\bZ}-{\rm characteristic}\}=1+\max\{ \lan S_{p,1}^{-1}x,x\ran\,
;\;\; x\in 2{\bZ}+1\}
\]\
\[
1-\min\{ \frac{x^2}{p};\;\; x\in 2{\bZ}+1 \}
=1-\frac{1}{p}\stackrel{(\ref{eq: liviu0})}{=} {\bf Froy}(L(p,1)).
\]
(b) $L(p,2)$  Again we distinguish two cases, $p=4k-1$ and
$p=4k+1$. Observe that
\[
S_{4k-1,2}=\left[
\begin{array}{cc}
-2k & 1\\
1 & -2
\end{array}
\right]
\]
is even so that $ \Theta_{4k-1,2}={\rm
rank}\,(S_{4k-1,2})=2\stackrel{(\ref{eq: liviu2})}{=}{\bf
Froy}(L(p,q))$.

When $p=4k+1$ we have
\[
S_{4k+1,2}=\left[
\begin{array}{cc}
-(2k+1) & 1 \\ 1 & -2
\end{array}
\right]
\]
with inverse
\[
S^{-1}_{4k+1,2}=\frac{1}{4k+1}\left[
\begin{array}{cc}
-2 & -1 \\
-1 & -(2k+1)
\end{array}
\right].
\]
The characteristic vectors  of  $S_{4k+1,2}$ are   determined by
the congruence
\[
\vec{v}=\left[
\begin{array}{c}
v^1\\
v^2
\end{array}
\right]\equiv \left[
\begin{array}{c}
2k+1\\
2
\end{array}
\right]\;\;\mod \;2.
\]
In particular, the vector
\[
\vec{u}_0=\left[
\begin{array}{c}
1\\
0
\end{array}
\right]
\]
is characteristic. We deduce
\[
{\bf Froy}(L(4k+1,2))\geq \Theta_{4k+1,2}=2+\max\{ \lan
S^{-1}_{4k+1,2}\vec{v}, \vec{v}\ran;\;\; \vec{v}-{\rm
characteristic}\}
\]
\[
\geq 2 + \lan S^{-1}_{4k+1,2}\vec{u}_0,\vec{u}_0\ran =
2-\frac{2}{4k+1} \stackrel{(\ref{eq: liviu2})}{=}{\bf
Froy}(L(4k+1,2)).
\]
The proposition is proved. $\Box$

\bigskip

 The Hirzebr\"{u}ch resolution is not minimal.
It can be transformed into a minimal one by  blowing-down
$-1$-spheres.  This operation changes the intersection form by a
diagonal, unimodular intersection  form and thus leaves the Elkies
invariant unchanged. Hence in the Statement of {\bf Conjecture 2} we can replace  $\Theta_{p,q}$ with $\Theta(S_{p,q}^{min})$ where $S_{p,q}^{min}$ denotes the intersection form of the minimal resolution.

The phenomenon claimed in the above conjecture and illustrates in
Proposition \ref{prop: conj2} is not singular. It was also
remarked in \cite{N} for  a large class of Brieskorn spheres. We
venture to formulate the following more general  statement.

\bigskip

\noindent {\bf Conjecture 3} {\em Suppose the rational homology sphere
$N$ is the link of an isolated complex singularity.  Denote   by
$q_{min}$ the intersection form  of the minimal resolution.  Then
for  any negative definite smooth manifold $X$ bounding $N$  we
have
\[
\Theta(q_X) \leq \Theta (q_{min}) = {\bf Froy}(N).
\]
Thus, loosely speaking, the minimal resolution is the most
complicated  smooth $4$-manifold bounding $N$.}

\section{The Casson-Walker invariant}
\setcounter{equation}{0}
 In this   section we describe  a
relationship between the Seiberg-Witten invariants of a  lens
space and its Casson-Walker invariant.

\subsection{The Seiberg-Witten invariants of a  rational homology sphere}
We use the same notations and terminology as in \S 2.1. Suppose $N$ is a rational
homology. The  set  $Spin^c(N)$ of s$spin^c$ structures on $N$ is
finite and has the same cardinality as $H:=H_1(N, {\bZ}$.  Fix a
$spin^c$ structure  $\sigma$ on $N$ and a good metric $g$.  Then
the set  of gauge equivalence classes of monopoles is finite. It
consists of an unique nondegenerate reducible monopole $\theta= (0, A_\si)$ and
finitely many,  nondegenerate  irreducible  ones $\{
\co_k;\;\;i=1,\ldots, n\}$.   Set
\[
n_k = i(\co_k, \theta)
\]
and $F(\si)= 4\eta(\dir_{A_\si})+\eta_{sign}$. The Seiberg-Witten
invariant of  $(N,\si)$ is the rational number
\be
{\bf sw}(\si)=\frac{1}{8} F(\si)-\sum_k (-1)^{n_k} .
\label{eq: def1}
\ee
In \cite{Chen1} and \cite{Lim} it  was proved that ${\bf sw}(\si)$
is independent of the choice of the good metric $g$ and
\[
{\bf sw}(\si) \in \frac{1}{8h_1} {\bZ}
\]
where $h_1=\# H_1(N, {\bZ})$. Observe that ${\bf sw}(\si)={\bf
sw}(\bar{\si})$ where $\si\mapsto \bar{\si}$ is the natural
involution on  $Spin^c(N)$.  Set
\be
{\bf sw}(N):=\sum_\si {\bf sw}(\si).
 \label{eq: def2}
\ee
If $N$ is a lens space $L(p,q)$ then, as explained  in \S 2.2,
a geometric Seifert structure on $N$ determines a
$spin^c$-structure $\si_0$ on $N$.    Will work with the geometric
Seifert structure determined by $\alpha=p/g.c.d.(p,q-1)$ and we set
\[
{\bf SW}_{p,q}=\sum_{j=0}^{p-1}{\bf sw}(\si_0\cdot t)t^j
\]
where $t$  is a generator of the cyclic group  ${\bZ}_p$. Observe that
\[
{\bf sw}(L(p,q))={\bf aug}({\bf SW}_{p,q}).
\]

The {\em Casson-Walker} invariant   of  $N$ is defined in
\cite{Lescop} and \cite{Walker}. It is a rational number $CW(N)$
uniquely determined by   certain Dehn surgery properties.

 We will work with  C. Lescop's normalization  used in
\cite{Lescop}.  It is related to    K. Walker's normalization used in
\cite{Walker} by the equality (\cite[Property T5.0, p.76]{Lescop}
\[
CW(N)_{Lescop}= \frac{h_1}{2}CW(N)_{Walker}.
\]
\begin{remark}{\rm The Casson-Walker invariant with C. Lescop's
normalization   differs by a sign  from the   conventions for the
Casson invariant used in \cite{FinSt} and \cite{N2}.  In these
references  the Casson invariant is normalized so that for the
Brieskorn homology sphere $\Sigma(a,b,c)$ we have
\[
CW(\Sigma(a,b,c)) = \frac{1}{8}\si(a,b,c)
\]
where $\si(a,b,c)$ denotes the signature of the Milnor fiber
associated to $\Sigma(a,b,c)$. In particular for the Poincar\'{e}
sphere $\Sigma(2,3,5)$ we have $\si(2,3,5 )={\rm sign}\,(-E_8)=-8$
so that the above formula gives the value $-1$ for the  Casson
invariant.     On the other hand  using C. Lescop's formula  \cite[Prop. 6.1.1,
To.0, p.97]{Lescop}  we obtain the value $1$ for the Casson-Walker
invariant. This explains the sign difference between the
definition of ${\bf sw}$ in (\ref{eq: def1}) and (\ref{eq: def2}) and the definition in
\cite{N2}.}
\end{remark}

The  Casson-Walker invariant of the lens space  can be expressed
in terms of the Dedekind sums. More precisely we have the equality
(see \cite{BL}, \cite{Walker})
\be
CW(L(p,q))=-\frac{p}{2}s(q,p).
\label{eq: walker}
\ee
We can now state the main result of this section.
\begin{theorem}
\[
{\bf sw}(L(p,q))=CW(L(p,q)).
\]
\label{th: cw}
\end{theorem}

\subsection{Seiberg-Witten $\Rightarrow$ Casson-Walker}
Our proof of Theorem \ref{th: cw} is arithmetic in nature and
relies on the computations in \S 2.2.

We will work with the same metric as in \S 2.2 and, since it has
positive scalar curvature we deduce there are no irreducible
monopoles, the unique reducible is also nondegenerate and thus
\[
{\bf sw}(L(p,q), \si) =\frac{1}{8} F_{p,q}(\si),\;\;\forall \si
\in Spin^c(L(p,q)).
\]
To proceed further we need to organize the   computational facts
established in \S 2.2 in a form suitable to  our current purposes.

Set $n= g.c.d.(p,q-1)$, $\alpha=p/n$
\[
\beta_2\cdot \frac{q-1}{n}\equiv -1 \;\mod\; \alpha,\;\;\beta_1=n-\beta_2
\]
\[
\omega_i=-\beta_i, \; \;q_i\omega_i\equiv 1\;\mod\;\alpha\;\;\forall\; i=1,2.
\]
The rational Euler degree of $L(p,q)$  equipped with the above
geometric  Seifert structure is
\[
\ell =-\frac{n}{\alpha}=-\frac{n^2}{p}.
\]
For each positive integer $m$ set $I_m:=\{0,1,\ldots, m-1\}$ and
$I^*_m=\{1,\ldots, m-1\}$. The set $Spin^c(L(p,q))$ can be
identified with $I_n\times I_\alpha$ and we have several functions
of interest
\[
\rho: I_n\times I_\alpha\ra {\bQ},\;\;\rho(k,j)=\frac{k}{n}
\]
\[
\gamma_1,\gamma_2:I_n\times I_\alpha\ra
{\bZ},\;\;\gamma_1(k,j)=j,\;\;\gamma_2(k,j)=k-1-j.
\]
The function $F_{p,q}(\si)$ can be regarded as a function
$F:I_n\times I_\alpha\ra {\bQ}$. It is explicitly described by
\[
F(k,j)=\ell +1 -4\ell\rho(1-\rho)+8\rho
-4\sum_{i=1}^2s(\omega_i,\alpha)-8\sum_{i=1}^2s(\omega_i,\alpha,
\frac{\gamma_i+\omega_i\rho}{\alpha}, -\rho)
\]
\be
+4\left\{
\begin{array}{rcl}
-\sum_{i=1}^2\LL\frac{q_i\gamma_i}{\alpha}\RR & {\rm if} & \rho =0
\\
 & & \\

(1-\frac{1}{\alpha})(1-2\rho)-\sum_{i=1}^2\{\frac{q_i\gamma_i+\rho}{\alpha}\}&
{\rm if} & \rho\neq 0
\end{array}
\right.
\label{eq: dedekind}
\ee
We have to prove
\be
\sum_{k\in I_n}\sum_{j\in I_\alpha} F(k,j) =-4ps(q,p).
\label{eq: cw1}
\ee
The proof of (\ref{eq: cw1}) relies on two identities. The first
one   was proved by M. Ouyang, \cite[p.652]{O}. More precisely, we
have
\be
\sum_{i=1}^2
s(\omega_i,\alpha)=s(q,p)-\frac{1}{6p}-\frac{n^2}{12p}+\frac{1}{4}.
\label{eq: ou}
\ee
The second one  is   central in the theory of  Dedekind sums  and
has the form
\be
\sum_{\mu\in I_m}\LL\frac{\mu+w}{m}\RR =((w)),\;\;
\forall\; m\in {\bZ}_+,\; w\in {\bR}.
\label{eq: kub}
\ee
For a proof we refer to \cite{HZ}.

Summing (\ref{eq: ou}) over $(k,j)\in I_n\times I_\alpha$ and using the equality $p=n\alpha$ we deduce
\be
4\sum_{k\in I_n}\sum_{j\in
I_\alpha}s(\omega_i,\alpha)=4ps(q,p)-\frac{2}{3}-\frac{n^2}{3}+p.
\label{eq: cw2}
\ee
We now proceed to sum over $(k,j)\in I_n\times I_\alpha$ all the
terms entering into the definition of $F(k,j)$.
\be
\sum_{k\in I_n}\sum_{j\in I_\alpha}(\ell+ 1)= -n^2+p.
\label{eq: cw3}
\ee
\be
8\sum_{k\in I_n}\sum_{j\in I_\alpha}\rho=8\sum_{j\in
I_\alpha}\sum_{k\in
I_n}\frac{k}{n}=\frac{8\alpha}{n}\frac{n(n-1}{2}=4(p-\alpha).
\label{eq: cw4}
\ee
\[
4\ell\sum_{k\in I_n}\sum_{j\in
I_\alpha}\rho(1-\rho)=-\frac{4n}{\alpha}\sum_{j\in I_\alpha}\sum_{k\in
I_n}\frac{k(n-k)}{n^2}=-\frac{4}{n}\sum_{k\in I_n}k(n-k)
\]
($\sum_{k\in I_n}k^2=\frac{n^3}{3}-\frac{n^2}{2}+\frac{n}{6}$)
\be
=-\frac{4}{n}(\frac{n^3}{2}-\frac{n^2}{2}-\frac{n^3}{3}+\frac{n^2}{2}-\frac{n}{6})=-\frac{2}{3}n^2+\frac{2}{3}.
\label{eq: cw5}
\ee
\[
\sum_{k\in I_n}\sum_{j\in I_\alpha}
S(\omega_i,\alpha,\frac{\gamma_i+\omega_i\rho}{\alpha},-\rho)=\sum_{k\in
I_n}\sum_{j\in I_\alpha}\sum_{\mu\in
I_\alpha}\LL\frac{\mu-\rho}{\alpha}\RR\LL\frac{(\omega_i(\mu-\rho)+\gamma_i+\omega_i\rho}{\alpha}\RR
\]
\[
\sum_{\mu\in I_\alpha}\LL\frac{\mu-\rho}{\alpha}\RR\sum_{k\in
I_n}\sum_{j\in
I_\alpha}\LL\frac{\gamma_i(k,j)+\omega_i\mu}{\alpha}\RR=\sum_{\mu\in I_\alpha}\LL\frac{\mu-\rho}{\alpha}\RR\sum_{k\in
I_n}\sum_{r\in
I_\alpha}\LL\frac{r+\omega_i\mu}{\alpha}\RR.
\]
According to (\ref{eq: kub}), the last sum (over $r$) is equal to
$((\omega_i\mu))=0$. Hence
\be
\sum_{k\in I_n}\sum_{j\in I_\alpha}
S(\omega_i,\alpha,\frac{\gamma_i+\omega_i\rho}{\alpha},-\rho)=0.
\label{eq: cw6}
\ee
Using (\ref{eq: kub}) again we deduce
\be
\sum_{k\in I_n^*}\sum_{j\in
I_\alpha}\LL\frac{q_i\gamma_i(k,j)}{\alpha}\RR=\sum_{k\in
I_n^*}\sum_{r\in
I_\alpha}\LL\frac{r}{\alpha}\RR\stackrel{(\ref{eq: kub})}{=}0.
\label{eq: cw7}
\ee
Observe that since $1-2\rho(k)=-(1-2\rho(n-k))$ we have
\be
(1-\frac{1}{\alpha})\sum_{k\in I_n^*}\sum_{j\in I_\alpha}(1-2\rho)=0.
\label{eq: cw8}
\ee
Finally  we have
\[
\sum_{k\in I_n^*}\sum_{j\in
I_\alpha}\left\{\frac{q_i\gamma_i(k,j)+\rho(k)}{\alpha}\right\}=\sum_{k\in
I_n}\sum_{j\in
I_\alpha}\left\{\frac{q_i\gamma_i(k,j)+\rho(k)}{\alpha}\right\}-\sum_{j\in
I_\alpha}\left\{\frac{q_i\gamma_i(0,j)}{\alpha}\right\}
\]
\[
=\sum_{k\in I_n}\sum_{j\in I_\alpha}\left\{
 \frac{nq_i\gamma_i(k,j) + k}{p}\right\}-\sum_{j\in
 I_\alpha}\left\{\frac{q_i\gamma_i(0,j)}{\alpha}\right\}.
 \]
 Now observe that as $k$ covers $I_n$ and $j$  covers $I_\alpha$ the
 quantity $(nq_i\gamma_i(k,j)+k\;\mod\;p)$   covers $I_p$ while
 $q_i\gamma_i(0,j)$  covers $I_\alpha$. Hence
 \[
\sum_{k\in I_n}\sum_{j\in I_\alpha}\left\{
 \frac{nq_i\gamma_i(k,j) + k}{p}\right\}=\sum_{r\in
 I_p}\{\frac{r}{p}\}= \frac{p-1}{2}
 \]
 and
\[
\sum_{j\in
 I_\alpha}\left\{\frac{q_i\gamma_i(0,j)}{\alpha}\right\}=\sum_{r\in
 I_\alpha} \{\frac{r}{\alpha}\} =\frac{\alpha-1}{2}.
 \]
 We conclude that
 \be
 4\sum_{k\in I_n^*}\sum_{j\in
 I_\alpha}\sum_{i=1}^{2}\left\{\frac{q_i\gamma_i(k,j)+\rho(k)}{\alpha}\right\}=
 4(p-\alpha).
 \label{eq: cw9}
 \ee
 The identity (\ref{eq: cw1}) now follows from (\ref{eq:
 cw2})-(\ref{eq: cw9}).   Theorem \ref{th: cw} is proved. $\Box$

 \bigskip

Taking into account Theorem \ref{th: cw} and the results of
\cite{Chen2} and \cite{Lim2} it is very tempting to formulate the
following conjecture.

\medskip

\noindent {\bf Conjecture 4} {\em For every rational homology
sphere  $N$  we have}
\[
{\bf sw}(N) =CW(N)
\]

\subsection{Seiberg-Witten invariants and Milnor torsion}
Consider the  invariant $T_{p,q}$  of the lens space $L(p,q)$ described  in the introduction.    The goal of this section is to prove  the following result.

\begin{proposition}{\rm If $g.c.d.(p,q-1)=1$ then
\be
T_{p,q}(1-t)(1-t^q)\sim \one
\label{eq: turaev1}
\ee
i.e. $T_{p,q}\sim \tau_{p,q}$.}
\label{prop: turaev}
\end{proposition}

\noindent{\bf Proof} \hspace{.3cm}  For a  while we will not rely
on the assumption $g.c.d.(p,q-1)=1$. We will continue to use the
notations in the previous subsection. Thus $n=g.c.d.(p,q-1)$.

As explained in \S 2.2, each $(k,j)\in I_n\times I_\alpha\cong
I_{n,\alpha}$   defines a line bundle on $L_j$ on  $L(p,q)$ and
thus, via the first Chern class an element $e(k,j)=c_1(L_{k,j})\in
H^2(L(p,q), {\bZ} \cong {\bZ}_p$. Moreover, the correspondence
\[
e:I_n\times I_\alpha\ra {\bZ}_p,\;\; (k,j)\mapsto e(k,j)
\]
is a bijection.

\begin{lemma}{\rm There exists an isomorphism of abelian groups $H^2(L(p,q), {\bZ})\ra {\bZ}_p$ such that}
\[
e(k,j)= q(k-1)-(q-1)j\; \mod \; p.
\]
\label{lemma: chern1}
\end{lemma}

\noindent {\bf Proof of the lemma}\hspace{.3cm}  $H^2(L(p,q), {\bZ})$ is  torsion  so according to the results in $\S 1.1$  it can be described in terms of the chosen geometric Seifert structure as follows.

Consider  map  ${\bQ}\oplus {\bZ}_\alpha\oplus {\bZ}_\alpha\ra {\bQ}/{\bZ}$
\[
(d,\gamma_1,\gamma_2)\mapsto d-\frac{\gamma_1+\gamma_2}{\alpha}
\]
 and the element
 \[
 L_0= (-n,\omega_1,\omega_2) \in \ker \delta.
 \]
 Recall that $L_0$ describes a line $V$-bundle  over an  genus $0$ orbifold whose  associated circle bundle coincides with the lens space equipped with the chosen Seifert structure. Then
 \[
 H^2(L(p,q), {\bZ})\cong \ker \delta / {\bZ}[L_0].
 \]
 Now observe  that $\ker \delta/ {\bZ}[L_0]$ has the presentation
 \[
 0\ra {\bZ}^2\stackrel{A}{\ra}{\bZ}^2\ra  \ker \delta / {\bZ}[L_0]\ra 0
 \]
 where
 \[
 A=\left[
 \begin{array}{cc}
 -n & 0 \\
 \omega_1 & \alpha
 \end{array}
 \right].
 \]
 We let the reader verify  that
 \be
\left[
\begin{array}{cc}
1 & 0 \\
0 & p
\end{array}
\right] = \left[
\begin{array}{cc}
-1 & 1\\
q & 1-q
\end{array}
\right]\cdot A \cdot \left[
\begin{array}{cc}
y & -\alpha \\
-x & -\omega_2
\end{array}
\right]
\label{eq: diagonalize}
\ee
where
\[
y=-(q-1)/n\;\;{\rm and}\;\;  x=-\frac{\omega_2y+ 1}{\alpha}.
\]
This shows that indeed
\[
\ker \delta / {\bZ}[L_0]\cong {\bZ}_p.
\]
To each pair $(k,j)\in I_n\times I_\alpha$ it corresponds the line bundle $L_{k,j}$ with  Seifert data $(k-1,j,k-1-j)\in \ker\delta$. Its first Chern class  is the image of the vector $\vec{v}=(k-1,j)\in {\bZ}^2$ in the quotient ${\bZ}^2/A{\bZ}^2$. Using the equality (\ref{eq: diagonalize}) we deduce that this image is $(y_2\;\mod \; p)$ where
\[
\left[
\begin{array}{c}
y_1\\
y_2
\end{array}
\right] =\left[
\begin{array}{cc}
-1 & 1\\
q & 1-q
\end{array}
\right]\cdot \left[
\begin{array}{c}
k-1 \\
j
\end{array}
\right].
\]
This establishes the assertion in  the lemma. {\bf q.e.d.}

\bigskip

Denote by $c: {\bZ}_p\ra I_n\times I_\alpha$ the inverse of the map $e$ described in the above lemma.

\begin{lemma}{\rm We have the following equalities.

\noindent (i) If $n=1$ then $\alpha =p$ and
\[
c(m)= (0,-\omega_2m+\omega_1\;\mod\; p).
\]
\noindent (ii) If $n\geq 1$ then
\[
c(-1)=c(p-1)= (0,\alpha-1)\;\;{\rm and}\;\;c(-m)= c(p-m)=(r,( -m -s\omega_1)\;\;\mod\;\alpha),\;\;\forall m\in I_p
\]
where $r\in I_n$ and $s\in {\bZ}$ are such that $ns= (m-1)+r$ so that
\[
r=-(m-1)\;\;\mod\;n\;\;{\rm and}\;\; s=\left\lceil\frac{m-1}{n}\right\rceil
\]
where $\lceil x \rceil$ is the smallest integer  $\geq x$.}
\label{lemma: count}
\end{lemma}

\noindent {\bf Proof}\hspace{.3cm} We prove only part (i). The second part is left to the reader.

Observe that  when $n=1$ we have $I_n\times I_\alpha=\{0\}\times I_\alpha$.  Thus we can write $c(m)=(0,j)$ where
\[
m=-q-(q-1)j\;\mod \;p
\]
Since $\omega_2=(q-1)^{-1}\;\mod\; p$ we have the following $\mod \;p$ equalities
\[
\omega_2m=-q\omega_2-j=-(q-1 +1)\omega_2 -j=-\omega_2 -1 -j.
\]
The equality in (i) now follows form $\omega_1+\omega_2=-n=-1$. {\bf q.e.d}

\bigskip

Now we can write
\[
{\bf SW}_{p,q}=\frac{1}{8}\sum_{m\in I_p}F(c(m))t^m.
\]
Since $\Sigma\cdot (1-t)=0$ in ${\bQ}[{\bZ}_p]$  the equality (\ref{eq: turaev1}) is equivalent to
\[
{\bf SW}_{p,q}(1-t)(1-t^q)\sim \one.
\]
We will prove  a slightly stronger statement namely
\be
{\bf SW}_{p,q}(1-t)(1-t^q)=\one.
\label{eq: turaev2}
\ee

Let us introduce  the  polynomial
\[
f(t)=\sum_{j\in I_p}\LL\frac{j}{p}\RR t^j\in {\bQ}[{\bZ}_p].
\]
A simple computation shows that
\[
f(t^{-1})=-f(t)
\]
and for all $m$ coprime with $p$ we have
\be
\left(\,\frac{1}{2}-f(t^m)\,\right)(1-t^m)= \one\;\;{\rm in}\;\;{\bQ}[{\bZ}_p]
\label{eq: inverse}
\ee
We want to express ${\bf SW}_{p,q}$ as a linear combinations of polynomials of the form $t^af(t^a)$, $t^af(t^a)f(t^b)$ and $\Sigma$. Observe first that since $n=1$,  in the equality (\ref{eq: dedekind}) of \S 3.2 we always have  $\rho =0$. Thus for all $(k,j)\in I_n\times I_\alpha$ we have
\[
F(k,j)=\ell + 1 -4\sum_{i=1}^2s(\omega_i, \alpha)
\]
\[
-8\sum_{i=1}^2s(\omega_i, \alpha,\gamma_i(k,j)/\alpha,0)-4\sum_{i=1}^2\LL\frac{q_i\gamma_i(k,j)}{\alpha}\RR.
\]
Observe two things.

\medskip

\noindent $\bullet$   Since $n=1$ we always have $k=0\in I_1=\{0\}$ so that we can write $\gamma_1(j)$ instead of $\gamma_i(k,j)$.

\noindent $\bullet$ The first term  in the definition of $F(k,j)$ is independent of $(k,j)$. Thus its contribution to  ${\bf SW}_{p,q}$  will be of the form $const.\Sigma$  which is cancelled upon multiplication  by $(1-t)$. Thus when computing ${\bf SW}_{p,q}(1-t)(1-t^q)$  we can neglect this first term.

\medskip

For $i=1,2$ define
\[
A_i=-8\sum_{m\in I_p}s\left(\omega_i, \alpha, \frac{\gamma_i(c(m))}{\alpha},0\right))t^m,\;\;B_i=\sum_{m\in
I_p}\LL\frac{q_i\gamma_i(c(m))}{\alpha}\RR t^m
\]
where according to \S 3.2 we have
\[
\gamma_1(j)=j,\;\;\gamma_2(j)=-1-j
\]
so that according to Lemma \ref{lemma: count} we have
\[
\gamma_1(c(m))= -\omega_2m+\omega_1,\;\;\gamma_2(c(m))=\omega_2m-\omega_1-1=\omega_2(m+1).
\]
Observe that  since $q_2\omega_2 =1\;\mod \; p$  and $\omega_2(q-1)=1\;\mod \; p$ we have
\[
q_2=(q-1)\;\mod\;p.
\]
\begin{lemma}
\be
B_1=-t^{-q}f(t^{-q})
\label{eq: b1}
\ee
\be
B_2=-t^{-1}f(t^{-1})
\label{eq: b2}
\ee
\be
A_1=-t^{-q}f(t^{-q})f(t^{q-1})
\label{eq: a1}
\ee
\be
A_2=t^{-1}f(t^{-1})f(t^{q-1})
\label{eq: a2}
\ee
\label{lemma: ab}
\end{lemma}

\noindent {\bf Proof} For any $(m,p+1$ we will denote by $1/m$ the inverse  of $m$ mod $p$.
\[
B_1=-\sum_{m\in I_m}\LL\frac{q_1(\omega_2m-\omega_1)}{\p}\RR t^m
\]
($\mu:=q_1\omega_2-q_1\omega_1=q_1\omega_2m-1$, $m=\frac{\omega_1}{\omega_2}(\mu+1)$ )
\[
=-t^{\omega_1/\omega_2}\sum_{\mu\in I_p}\LL\frac{\mu}{p}\RR t^{\omega_1\mu/\omega_2} = -t^{\omega_1/\omega_2}f(t^{\omega_1/\omega_2}).
\]
Now  observe that $1/\omega_2=q_2=q-1$ and $\omega_1=-1-\omega_2$ so that $\omega_1/\omega_2 =-q$. This proves (\ref{eq: b1}).
\[
B_2=\sum_{m\in I_m}\LL\frac{q_2\omega_2(m+1)}{p}\RR t^m=\sum_{\mu\in I_p}\LL\frac{\mu}{p}\RR
t^{\mu-1} = t^{-1}f(t)=-t^{-1}f(t^{-1}).
\]
This proves (\ref{eq: b2}).
\[
A_1=\sum_{m\in I_p}\sum_{\mu\in I_p}\LL\frac{\mu}{p}\RR\LL\frac{\omega_1\mu-\omega_2m+\omega_1}{p}\RR t^m =\sum_{\mu\in I_p}\LL\frac{\mu}{p}\RR\sum_{m\in I_p}\LL\frac{\omega_1\mu-\omega_2m+\omega_1}{p}\RR t^m
\]
( $r=\omega_1\mu-\omega_2m+\omega_1$, $m=-r/\omega_2+\omega_1(\mu+1)/\omega_2$ )
\[
=t^{\omega_1/\omega_2}\sum_{\mu\in I_p}\LL\frac{\mu}{p}\RR t^{\omega_1\mu/\omega_2}\sum_{r\in I_p}\LL\frac{r}{p}\RR t^{-r/\omega_2}=t^{\omega_1/\omega_2}f(t^{\omega^1/\omega_2})f(t^{-1/\omega_2})
\]
\[
=t^{-q}f(t^{-q})f(t^{-(q-1)})=-t^{-q}f(t^{-q})f(t^{q-1}).
\]
This proves (\ref{eq: a1}). Finally, we have
\[
A_2=\sum_{m\in I_p}\sum_{\mu\in I_p}\LL\frac{\mu}{p}\RR \LL\frac{\omega_2\mu+\omega_2m+\omega_2}{p}\RR t^m = \sum_{\mu\in I_p}\LL\frac{\mu}{p}\RR\sum_{m\in I_p}\LL\frac{\omega_2\mu+\omega_2m+\omega_2}{p}\RR t^m
\]
( $r=\omega_2(m+\mu+1)$, $m=r/\omega_2-\mu-1$ )
\[
=t^{-1}\sum_{\mu\in I_p}\LL\frac{\mu}{p}\RR t^{-\mu}\sum_{r\in I_p}\LL\frac{r}{p}\RR t^{r/\omega_2}=t^{-1}f(t^{-1})f(t^{q-1})
\]
This proves (\ref{eq: a2}). {\bf q.e.d.}

\bigskip

We can now finish the proof of Proposition \ref{prop: turaev}. Using Lemma \ref{lemma: ab}  we deduce
\[
8{\bf SW}_{p,q}(1-t)(1-t^q)=(-8A_1-8A_2-4B_1-B_2 +const.\Sigma)(1-t)(1-t^q)
\]
\[
=-4(2A_1+2A_2+B_1+B_2)(1-t)(1-t^q)
\]
\[
=-4\left\{\;- t^{-q}f(t^{-q})(\,1+2f(t^{q-1})\;)-t^{-1}f(t^{-1})(\, 1-2f(t^{q-1}\,)\,)\;\right\}(1-t)(1-t^q)
\]
\[
=-8\left\{\;-t^{-q}f(t^{-q})(\,\frac{1}{2}-f(t^{-(q-1)})\,)-t^{-1}f(t^{-1})(\frac{1}{2}-f(t^{q-1})\,)\;\right\}(1-t)(1-t^q)
\]
\[
\stackrel{(\ref{eq:
inverse})}{=}8\left\{ t^{-q}f(t^{-q})\cdot\frac{\one}{1-t^{1-q}}+t^{-1}f(t^{-1})\cdot\frac{\one}{1-t^{q-1}}\;\right\}(1-t)(1-t^q)
\]
\[
=8\left\{ t^{-1}f(t^{-q})\cdot\frac{\one}{t^{q-1}-1}+t^{-1}f(t^{-1})\cdot\frac{\one}{1-t^{q-1}}\;\right\}(1-t)(1-t^q)
\]
\[
=8t^{-1}\frac{\one}{1-t^{q-1}}(f(t^{-1})-f(t^{-q})\,)(1-t)(1-t^q)
\]
\[
=8t^{-1}\frac{\one}{1-t^{q-1}}(f(t^q)-f(t)\,)(1-t)(1-t^q)
\]
\[
\stackrel{(\ref{eq: inverse})}{=}8t^{-1}\frac{\one}{1-t^{q-1}}\left(\frac{\one}{1-t}-\frac{\one}{1-t^q}\right)(1-t)(1-t^q)
\]
\[
=8t^{-1}\frac{\one}{1-t^{q-1}}\left\{\; (1-t^q)-(1-t)\;\right\} = 8t^{-1}\frac{\one}{1-t^{q-1}}(t-t^q)= 8\cdot \one.
\]
The proof of Proposition \ref{prop: turaev} is now complete. $\Box$

\bigskip

The restriction $g.c.d.(p,q-1)=1$  in Proposition \ref{prop:
turaev}  can be dropped but we will present the details elsewhere.
Instead, we  have included  below    an explicit description of $T_{p,q}$ for all
$1<q<p\leq 10$. These computations    confirm the validity of
(\ref{eq: turaev1}) even if $g.c.d.(p,q-1)>1$.

Proposition \ref{prop: turaev} confirms  a hypothesis formulated in  \cite{Turaev2}.    We formulate the
following  conjecture.

\bigskip

\noindent {\bf Conjecture 5} {\em For  any rational homology
sphere $N$  the  augmentation-free part of the Seiberg-Witten
invariant coincides with the refined torsion  of Turaev,
\cite{Turaev2}.}

\bigskip

\begin{center}
{\bf Numerical experiments}
\end{center}

Below we let the reader verify  the elementary  identity (\ref{eq:
turaev1}) in each case.

\noindent $\bullet$ L(2,q)

\[
{T_{2, \,1}} \sim  - {\displaystyle \frac {1}{8}} \,t +
{\displaystyle \frac {1}{8}}
\]

\noindent $\bullet$ L(3,q)

\[
{T_{3, \,1}} \sim {\displaystyle \frac {1}{9}} \,t^{2} -
{\displaystyle \frac {2}{9}} \,t + {\displaystyle \frac {1}{9}}
\]

\[
{T_{3, \,2}} \sim  - {\displaystyle \frac {1}{9}} \,t^{2} +
{\displaystyle \frac {2}{9}} \,t - {\displaystyle \frac {1}{9}}
\]

\noindent $\bullet$ L(4,q)

\[
{T_{4, \,1}} \sim  - {\displaystyle \frac {5}{16}} \,t^{3} +
{\displaystyle \frac {1}{16}} \,t^{2} + {\displaystyle \frac {3}{
16}} \,t + {\displaystyle \frac {1}{16}}
\]

\[
{T_{4, \,3}} \sim  - {\displaystyle \frac {5}{16}} \,t^{3} +
{\displaystyle \frac {1}{16}} \,t^{2} + {\displaystyle \frac {3}{
16}} \,t + {\displaystyle \frac {1}{16}}
\]

\noindent $\bullet$ L(5,q)

\[
{T_{5, \,1}} \sim {\displaystyle \frac {1}{5}} \,t^{4} -
{\displaystyle \frac {2}{5}} \,t^{2} + {\displaystyle \frac {1}{5
}}
\]

\[
{T_{5, \,2}} \sim  - {\displaystyle \frac {1}{5}} \,t^{4} +
{\displaystyle \frac {1}{5}} \,t^{3} + {\displaystyle \frac {1}{5
}} \,t - {\displaystyle \frac {1}{5}}
\]

\[
T_{5,\,3}\sim
-\frac{1}{5}t^4+\frac{1}{5}t^3+\frac{1}{5}t-\frac{1}{5}
\]

\[
T_{5,\,4}\sim -\frac{1}{5}t^3 +\frac{2}{5}t^2 -\frac{1}{5}t
\]

\noindent $\bullet$ L(6,q)

\[
{T_{6, \,1}} \sim  - {\displaystyle \frac {35}{72}} \,t^{5} -
{\displaystyle \frac {5}{72}} \,t^{4} + {\displaystyle \frac {13
}{72}} \,t^{3} + {\displaystyle \frac {19}{72}} \,t^{2} +
{\displaystyle \frac {13}{72}} \,t - {\displaystyle \frac {5}{72}
}
\]

\[
{T_{6, \,5}} \sim {\displaystyle \frac {35}{72}} \,t^{5} +
{\displaystyle \frac {5}{72}} \,t^{4} - {\displaystyle \frac {13
}{72}} \,t^{3} - {\displaystyle \frac {19}{72}} \,t^{2} -
{\displaystyle \frac {13}{72}} \,t + {\displaystyle \frac {5}{72}
}
\]

\noindent $\bullet$ L(7,q)

\[
{T_{7, \,1}} \sim {\displaystyle \frac {2}{7}} \,t^{6} +
{\displaystyle \frac {1}{7}} \,t^{5} - {\displaystyle \frac {1}{7
}} \,t^{4} - {\displaystyle \frac {4}{7}} \,t^{3} -
{\displaystyle \frac {1}{7}} \,t^{2} + {\displaystyle \frac {1}{7
}} \,t + {\displaystyle \frac {2}{7}}
\]

\[
{T_{7, \,2}} \sim  - {\displaystyle \frac {2}{7}} \,t^{6} +
{\displaystyle \frac {1}{7}} \,t^{5} + {\displaystyle \frac {2}{7
}} \,t^{3} + {\displaystyle \frac {1}{7}} \,t - {\displaystyle
\frac {2}{7}}
\]

\[
{T_{7, \,3}} \sim  - {\displaystyle \frac {1}{7}} \,t^{6} +
{\displaystyle \frac {2}{7}} \,t^{4} - {\displaystyle \frac {2}{7
}} \,t^{3} + {\displaystyle \frac {2}{7}} \,t^{2} -
{\displaystyle \frac {1}{7}}
\]

\[
{T_{7, \,4}} \sim  - {\displaystyle \frac {2}{7}} \,t^{6} +
{\displaystyle \frac {1}{7}} \,t^{5} + {\displaystyle \frac {2}{7
}} \,t^{3} + {\displaystyle \frac {1}{7}} \,t - {\displaystyle
\frac {2}{7}}
\]

\[
{T_{7, \,5}} \sim  - {\displaystyle \frac {1}{7}} \,t^{6} +
{\displaystyle \frac {2}{7}} \,t^{4} - {\displaystyle \frac {2}{7
}} \,t^{3} + {\displaystyle \frac {2}{7}} \,t^{2} -
{\displaystyle \frac {1}{7}}
\]

\[
{T_{7, \,6}} \sim {\displaystyle \frac {1}{7}} \,t^{6} -
{\displaystyle \frac {2}{7}} \,t^{5} - {\displaystyle \frac {1}{7
}} \,t^{4} + {\displaystyle \frac {4}{7}} \,t^{3} -
{\displaystyle \frac {1}{7}} \,t^{2} - {\displaystyle \frac {2}{7
}} \,t + {\displaystyle \frac {1}{7}}
\]

\noindent $\bullet$ L(8,q)

\[
{T_{8, \,1}} \sim  - {\displaystyle \frac {21}{32}} \,t^{7} -
{\displaystyle \frac {7}{32}} \,t^{6} + {\displaystyle \frac {3}{
32}} \,t^{5} + {\displaystyle \frac {9}{32}} \,t^{4} +
{\displaystyle \frac {11}{32}} \,t^{3} + {\displaystyle \frac {9
}{32}} \,t^{2} + {\displaystyle \frac {3}{32}} \,t -
{\displaystyle \frac {7}{32}}
\]

\[
{T_{8, \,3}} \sim  - {\displaystyle \frac {9}{32}} \,t^{7} -
{\displaystyle \frac {3}{32}} \,t^{6} - {\displaystyle \frac {9}{
32}} \,t^{5} + {\displaystyle \frac {5}{32}} \,t^{4} +
{\displaystyle \frac {7}{32}} \,t^{3} - {\displaystyle \frac {3}{
32}} \,t^{2} + {\displaystyle \frac {7}{32}} \,t +
{\displaystyle \frac {5}{32}}
\]

\[
{T_{8, \,5}} \sim  - {\displaystyle \frac {5}{32}} \,t^{7} -
{\displaystyle \frac {7}{32}} \,t^{6} + {\displaystyle \frac {3}{
32}} \,t^{5} - {\displaystyle \frac {7}{32}} \,t^{4} -
{\displaystyle \frac {5}{32}} \,t^{3} + {\displaystyle \frac {9}{
32}} \,t^{2} + {\displaystyle \frac {3}{32}} \,t +
{\displaystyle \frac {9}{32}}
\]

\[
{T_{8, \,7}} \sim {\displaystyle \frac {21}{32}} \,t^{7} -
{\displaystyle \frac {9}{32}} \,t^{6} - {\displaystyle \frac {3}{
32}} \,t^{5} + {\displaystyle \frac {7}{32}} \,t^{4} -
{\displaystyle \frac {11}{32}} \,t^{3} + {\displaystyle \frac {7
}{32}} \,t^{2} - {\displaystyle \frac {3}{32}} \,t -
{\displaystyle \frac {9}{32}}
\]

\noindent $\bullet$ L(9,q)

\[
{T_{9, \,1}} \sim {\displaystyle \frac {10}{27}} \,t^{8} +
{\displaystyle \frac {7}{27}} \,t^{7} + {\displaystyle \frac {1}{
27}} \,t^{6} - {\displaystyle \frac {8}{27}} \,t^{5} -
{\displaystyle \frac {20}{27}} \,t^{4} - {\displaystyle \frac {8
}{27}} \,t^{3} + {\displaystyle \frac {1}{27}} \,t^{2} +
{\displaystyle \frac {7}{27}} \,t + {\displaystyle \frac {10}{27}
}
\]

\[
{T_{9, \,2}} \sim  - {\displaystyle \frac {10}{27}} \,t^{8} +
{\displaystyle \frac {2}{27}} \,t^{7} - {\displaystyle \frac {1}{
27}} \,t^{6} + {\displaystyle \frac {8}{27}} \,t^{5} +
{\displaystyle \frac {2}{27}} \,t^{4} + {\displaystyle \frac {8}{
27}} \,t^{3} - {\displaystyle \frac {1}{27}} \,t^{2} +
{\displaystyle \frac {2}{27}} \,t - {\displaystyle \frac {10}{27}
}
\]

\[
{T_{9, \,4}} \sim {\displaystyle \frac {1}{27}} \,t^{8} -
{\displaystyle \frac {2}{27}} \,t^{7} - {\displaystyle \frac {8}{
27}} \,t^{6} + {\displaystyle \frac {10}{27}} \,t^{5} -
{\displaystyle \frac {2}{27}} \,t^{4} + {\displaystyle \frac {10
}{27}} \,t^{3} - {\displaystyle \frac {8}{27}} \,t^{2} -
{\displaystyle \frac {2}{27}} \,t + {\displaystyle \frac {1}{27}
}
\]

\[
{T_{9, \,5}} \sim  - {\displaystyle \frac {10}{27}} \,t^{8} +
{\displaystyle \frac {2}{27}} \,t^{7} - {\displaystyle \frac {1}{
27}} \,t^{6} + {\displaystyle \frac {8}{27}} \,t^{5} +
{\displaystyle \frac {2}{27}} \,t^{4} + {\displaystyle \frac {8}{
27}} \,t^{3} - {\displaystyle \frac {1}{27}} \,t^{2} +
{\displaystyle \frac {2}{27}} \,t - {\displaystyle \frac {10}{27}
}
\]

\[
{T_{9, \,7}} \sim  - {\displaystyle \frac {8}{27}} \,t^{8} -
{\displaystyle \frac {2}{27}} \,t^{7} + {\displaystyle \frac {10
}{27}} \,t^{6} + {\displaystyle \frac {1}{27}} \,t^{5} -
{\displaystyle \frac {2}{27}} \,t^{4} + {\displaystyle \frac {1}{
27}} \,t^{3} + {\displaystyle \frac {10}{27}} \,t^{2} -
{\displaystyle \frac {2}{27}} \,t - {\displaystyle \frac {8}{27}
}
\]

\[
{T_{9, \,8}} \sim {\displaystyle \frac {8}{27}} \,t^{8} -
{\displaystyle \frac {7}{27}} \,t^{7} - {\displaystyle \frac {10
}{27}} \,t^{6} - {\displaystyle \frac {1}{27}} \,t^{5} +
{\displaystyle \frac {20}{27}} \,t^{4} - {\displaystyle \frac {1
}{27}} \,t^{3} - {\displaystyle \frac {10}{27}} \,t^{2} -
{\displaystyle \frac {7}{27}} \,t + {\displaystyle \frac {8}{27}
}
\]

\noindent $\bullet$ L(10,q)

\[
{T_{10, \,1}} \sim  - {\displaystyle \frac {33}{40}} \,t^{9} -
{\displaystyle \frac {3}{8}} \,t^{8} - {\displaystyle \frac {1}{
40}} \,t^{7} + {\displaystyle \frac {9}{40}} \,t^{6} +
{\displaystyle \frac {3}{8}} \,t^{5} + {\displaystyle \frac {17}{
40}} \,t^{4} + {\displaystyle \frac {3}{8}} \,t^{3} +
{\displaystyle \frac {9}{40}} \,t^{2} - {\displaystyle \frac {1}{
40}} \,t - {\displaystyle \frac {3}{8}}
\]

\[
{T_{10, \,3}} \sim  - {\displaystyle \frac {3}{8}} \,t^{9} -
{\displaystyle \frac {1}{40}} \,t^{8} + {\displaystyle \frac {1}{
40}} \,t^{7} - {\displaystyle \frac {9}{40}} \,t^{6} +
{\displaystyle \frac {9}{40}} \,t^{5} + {\displaystyle \frac {3}{
8}} \,t^{4} + {\displaystyle \frac {9}{40}} \,t^{3} -
{\displaystyle \frac {9}{40}} \,t^{2} + {\displaystyle \frac {1}{
40}} \,t - {\displaystyle \frac {1}{40}}
\]

\[
{T_{10, \,7}} \sim  - {\displaystyle \frac {3}{8}} \,t^{9} -
{\displaystyle \frac {9}{40}} \,t^{8} + {\displaystyle \frac {9}{
40}} \,t^{7} - {\displaystyle \frac {1}{40}} \,t^{6} +
{\displaystyle \frac {1}{40}} \,t^{5} + {\displaystyle \frac {3}{
8}} \,t^{4} + {\displaystyle \frac {1}{40}} \,t^{3} -
{\displaystyle \frac {1}{40}} \,t^{2} + {\displaystyle \frac {9}{
40}} \,t - {\displaystyle \frac {9}{40}}
\]

\[
{T_{10, \,9}} \sim {\displaystyle \frac {33}{40}} \,t^{9} -
{\displaystyle \frac {9}{40}} \,t^{8} - {\displaystyle \frac {3}{
8}} \,t^{7} + {\displaystyle \frac {3}{8}} \,t^{6} +
{\displaystyle \frac {1}{40}} \,t^{5} - {\displaystyle \frac {17
}{40}} \,t^{4} + {\displaystyle \frac {1}{40}} \,t^{3} +
{\displaystyle \frac {3}{8}} \,t^{2} - {\displaystyle \frac {3}{8
}} \,t - {\displaystyle \frac {9}{40}}
\]


\begin{thebibliography}{XXXXX}

\addcontentsline{toc}{section}{References}

\bibitem{BPV} W. Barth, C. Peters, A. Van de Ven: {\sl  Compact
Complex Surfaces}, Erg. der Math., 2. Folge, Band 4, Springer
Verlag, Berlin, 1984.

\bibitem{BL} S. Boyer, D. Lines: {\sl Surgery formul{\ae} for
Casson's invariant and extensions to homology lens spaces}, J.
Reine. Angew. Math., {\bf 405}(1990), 181-220.

\bibitem{Chen1} W.Chen: {\sl Casson invariant and Seiberg-Witten
gauge theory}, Turkish J. Math., {\bf 21}(1997), 61-81.

\bibitem{Chen2} W. Chen: {\sl Dehn surgery formula for
Seiberg-Witten invariants of homology 3-spheres}, {\sf
dg-ga/9708006}

\bibitem{CH} A.J. Casson, J.L. Harer: {\sl Some homology lens spaces which bound rational homology balls},  Pacific J. Math., {\bf 96}(1981), 23-36.

\bibitem{Elk} N.D. Elkies: {\sl  A characterization of the ${\bZ}^n$ lattice}, Math. Res. Lett., {\bf 2}(1995), 321-326.

\bibitem{FinSt} R. Fintushel, R. Stern: {\sl  Instanton homology of
Seifert fibered homology three spheres}, Proc. London Math. Soc.
{\bf 61}(1990), 109-137.

\bibitem{Fr} K.A. Froyshov: {\sl The Seiberg-Witten equations and four manifolds with boundary}, Math. Res. Lett., {\bf 3} (1996), 373-390.

\bibitem{FS} M. Furuta, B. Steer: {\sl Seifert fibred homology 3-spheres and the Yang-Mills equations on Riemann surfaces with marked points},  Adv. in Math. {\bf 96}(1992),  38-102.


\bibitem{Gaulter} M. Gaulter: {\sl  Lattices without short characteristic vectors}, Math. Res. Lett., {\bf 5}(1998), 353-362.


\bibitem{HNK} F. Hirzebr\"{u}ch, W.D. Neumann, S.S. Koh: {\sl Differentiable Manifolds and Quadratic Forms}, Lect. Notes in Pure and Appl. Math., No. 4, Marcel Dekker, 1971.

\bibitem{HZ} F. Hirzebr\"{u}ch, D. Zagier: {\sl The Atiyah-Singer Index Theorem and Elementary Number Theory},  Math. Lect. Series {\bf 3}, Publish or Perish Inc., Boston,  1974.


\bibitem{JN} M. Jankins, W. D. Neumann: {\sl Lectures on Seifert Manifolds}, Brandeis Lecture Notes, 1983.


\bibitem{Lescop} C. Lescop: {\sl Global Surgery Formula for the
Casson-Walker Invariant}, Annals of Math. Studies, vol. {\bf 140}, Princeton University Press, 1996.
University Press

\bibitem{Lim} Y. Lim: {\sl  Seiberg-Witten invariants for 3-manifolds in the case $b_1=0$ or $1$}, preprint, 1998.

\bibitem{Lim2} Y. Lim: {\sl The equivalence of Seiberg-Witten and
Casson invariants for homology 3-spheres}, preprint.

\bibitem{Lisca} P. Lisca: {\sl Symplectic fillings and positive scalar curvature}, Geometry and Topology, {\bf 2}(1998), 103-116.


\bibitem{Mar} M. Marcolli: {\sl  Equivariant Seiberg-Witten-Floer
homology}, {\sf dg-ga 9606003}.

\bibitem{Milnor} J. Milnor: {\sl  Whitehead torsion},  Bull. Amer.
Math. Soc. {\bf 72}(1966), 358-426.

\bibitem{NR} W.D. Neumann, F. Raymond: {\sl Seifert manifolds, plumbing, $\mu$-invariant  and orientation reversing maps}, in {\sl``Lecture Notes in Mathematics''}, vol. {\bf 644}, 161-195.

\bibitem{N0} L.I. Nicolaescu: {\sl Adiabatic limits of the  Seiberg-Witten equations on Seifert manifolds},  Comm. Anal. and Geom., {\bf 6}(1998), 301-362.


\bibitem{N} L.I. Nicolaescu: {\sl  Finite energy Seiberg-Witten moduli spaces on $4$-manifolds bounding Seifert fibrations}, {\sf dg-ga 9711006}.

\bibitem{N1} L.I. Nicolaescu: {\sl Eta invariants of Dirac operators on  circle bundles over Riemann surfaces and virtual  dimensions of finite energy Seiberg-Witten moduli spaces}, {\sf math.DG/9805046}, Israel. J. Math, to appear.

\bibitem{N2} L.I. Nicolaescu: {\sl Lattice points, Dedekind-Rademacher sums and a conjecture of Kronheimer and
Mrowka},  {\sf math.DG/9801030}.

\bibitem{Or} P. Orlik: {\sl Seifert Manifolds}, Lect. Motes in Math., vol. {\bf 291}, Springer-Verlag, 1972.

\bibitem{O} M. Ouyang: {\sl  Geometric invariants for Seifert fibered 3-manifolds},  Trans. Amer. Math. Soc. {\bf 346}(1994), 641-659.

\bibitem{Ra} H. Rademacher: {\sl Some remarks on certain generalized Dedekind sums}, Acta Arithmetica, {\bf 9}(1964), 97-105.

\bibitem{RG} H. Rademacher, E. Grosswald: {\sl Dedekind Sums }, The Carus Math. Monographs, MAA, 1972.


\bibitem{Ran} R. von Randow: {\sl  Z\"{u}r Topologie von dreidimensionalen Baummanigfatigkeiten}, Bonner Math. Schriften, {\bf 14}(1962).


\bibitem{S} P. Scott: {\sl The geometries of 3-manifolds}, Bull. London. Math. Soc. {\bf 15}(1983), 401-487.

\bibitem{Turaev}  V.G. Turaev: {\sl  Euler structures, nonsingular
vector fields and torsions of Reidemeister type}, Izvestia Akad. Nauk. USSR, {\bf 53}(1989); English Transl.: Math. USSR
Izvestia, {\bf 34}(1990), 627-662

\bibitem{Turaev2} V.G. Turaev: {\sl  Torsion invariants of
$spin^c$ structures on 3-manifolds}, Math. Res. Letters, {\bf
4}(1997), 679-695.

\bibitem{Walker} K. Walker: {\sl  An Extension of Casson's
Invariant}, Annals of Math. Studies, vol. {\bf 126}, Princeton University
Press,  1992.


\end{thebibliography}
\end{document}